\documentclass[12pt]{amsart}

\usepackage{amsmath,amssymb,amsthm,amscd}

\begin{document}

\newcommand{\End}{{\rm{End}\ts}}
\newcommand{\Hom}{{\rm{Hom}}}
\newcommand{\ch}{{\rm{ch}\ts}}
\newcommand{\non}{\nonumber}
\newcommand{\wt}{\widetilde}
\newcommand{\wh}{\widehat}
\newcommand{\ot}{\otimes}
\newcommand{\la}{\lambda}
\newcommand{\La}{\Lambda}
\newcommand{\al}{\alpha}
\newcommand{\be}{\beta}
\newcommand{\ga}{\gamma}
\newcommand{\si}{\sigma}
\newcommand{\vp}{\varphi}
\newcommand{\de}{\delta^{}}
\newcommand{\om}{\omega^{}}
\newcommand{\hra}{\hookrightarrow}
\newcommand{\ve}{\varepsilon}
\newcommand{\ts}{\,}
\newcommand{\qin}{q^{-1}}
\newcommand{\tss}{\hspace{1pt}}
\newcommand{\U}{ {\rm U}}
\newcommand{\Y}{ {\rm Y}}
\newcommand{\AAb}{\mathbb{A}\tss}
\newcommand{\CC}{\mathbb{C}\tss}
\newcommand{\QQ}{\mathbb{Q}\tss}
\newcommand{\SSb}{\mathbb{S}\tss}
\newcommand{\ZZ}{\mathbb{Z}\tss}
\newcommand{\Z}{{\rm Z}}
\newcommand{\Ac}{\mathcal{A}}
\newcommand{\Pc}{\mathcal{P}}
\newcommand{\Qc}{\mathcal{Q}}
\newcommand{\Tc}{\mathcal{T}}
\newcommand{\Sc}{\mathcal{S}}
\newcommand{\Bc}{\mathcal{B}}
\newcommand{\Ec}{\mathcal{E}}
\newcommand{\Hc}{\mathcal{H}}
\newcommand{\Ar}{{\rm A}}
\newcommand{\Ir}{{\rm I}}
\newcommand{\Zr}{{\rm Z}}
\newcommand{\gl}{\mathfrak{gl}}
\newcommand{\Pf}{{\rm Pf}}
\newcommand{\oa}{\mathfrak{o}}
\newcommand{\spa}{\mathfrak{sp}}
\newcommand{\g}{\mathfrak{g}}
\newcommand{\ka}{\mathfrak{k}}
\newcommand{\p}{\mathfrak{p}}
\newcommand{\sll}{\mathfrak{sl}}
\newcommand{\agot}{\mathfrak{a}}
\newcommand{\qdet}{ {\rm qdet}\ts}
\newcommand{\sdet}{ {\rm sdet}\ts}
\newcommand{\Gr}{ {\rm Gr}\tss}
\newcommand{\sgn}{ {\rm sgn}\ts}
\newcommand{\Sym}{\mathfrak S}
\newcommand{\fand}{\quad\text{and}\quad}
\newcommand{\Fand}{\qquad\text{and}\qquad}
\newcommand{\vt}{{\tss|\hspace{-1.5pt}|\tss}}

\renewcommand{\theequation}{\arabic{section}.\arabic{equation}}

\newtheorem{thm}{Theorem}[section]
\newtheorem{lem}[thm]{Lemma}
\newtheorem{prop}[thm]{Proposition}
\newtheorem{cor}[thm]{Corollary}
\newtheorem{conj}[thm]{Conjecture}

\theoremstyle{definition}
\newtheorem{defin}[thm]{Definition}

\theoremstyle{remark}
\newtheorem{remark}[thm]{Remark}
\newtheorem{example}[thm]{Example}

\newcommand{\bth}{\begin{thm}}
\renewcommand{\eth}{\end{thm}}
\newcommand{\bpr}{\begin{prop}}
\newcommand{\epr}{\end{prop}}
\newcommand{\ble}{\begin{lem}}
\newcommand{\ele}{\end{lem}}
\newcommand{\bco}{\begin{cor}}
\newcommand{\eco}{\end{cor}}
\newcommand{\bde}{\begin{defin}}
\newcommand{\ede}{\end{defin}}
\newcommand{\bex}{\begin{example}}
\newcommand{\eex}{\end{example}}
\newcommand{\bre}{\begin{remark}}
\newcommand{\ere}{\end{remark}}
\newcommand{\bcj}{\begin{conj}}
\newcommand{\ecj}{\end{conj}}

\newcommand{\bal}{\begin{aligned}}
\newcommand{\eal}{\end{aligned}}
\newcommand{\beq}{\begin{equation}}
\newcommand{\eeq}{\end{equation}}
\newcommand{\ben}{\begin{equation*}}
\newcommand{\een}{\end{equation*}}

\newcommand{\bpf}{\begin{proof}}
\newcommand{\epf}{\end{proof}}

\def\beql#1{\begin{equation}\label{#1}}

\title[Comultiplication rules]
{Comultiplication rules for the double
Schur\\[0.2em] functions and Cauchy identities}

\author{A. I. Molev}

\address{School of Mathematics and Statistics\\
University of Sydney, NSW 2006, Australia}
\email{alexm\hspace{0.09em}@\hspace{0.1em}maths.usyd.edu.au}

\date{} 


\begin{abstract}
The double Schur functions form a distinguished basis
of the ring $\Lambda(x\vt a)$ which is a multiparameter
generalization of the ring of symmetric functions $\Lambda(x)$.
The canonical comultiplication on $\Lambda(x)$ is extended
to $\Lambda(x\vt a)$ in a natural way so that the double power
sums symmetric functions are primitive elements.
We calculate the dual Littlewood--Richardson coefficients
in two different ways thus
providing comultiplication rules for the
double Schur functions. We also
prove multiparameter analogues of the Cauchy identity.
A new family of Schur type functions plays the role
of a dual object in the identities.
We describe some properties of these dual
Schur functions including a combinatorial presentation
and an expansion formula in terms of the ordinary
Schur functions. The dual Littlewood--Richardson coefficients
provide a multiplication rule for the dual Schur functions.
\end{abstract}

\maketitle

\tableofcontents

\section{Introduction}
\label{sec:int}
\setcounter{equation}{0}

The ring $\La=\La(x)$ of symmetric functions in the
set of variables $x=(x_1,x_2,\dots)$ admits a multiparameter
generalization $\La(x\vt a)$, where $a$ is a sequence of variables
$a=(a_i)$, $i\in\ZZ$. Let $\QQ[a]$ denote the ring
of polynomials in the variables $a_i$ with rational coefficients.
The ring $\La(x\vt a)$ is generated over $\QQ[a]$ by
the {\it double power sums symmetric functions\/}
\beql{powess}
p_k(x\vt a)=\sum_{i=1}^{\infty} (x_i^k-a_i^k).
\eeq
Moreover, it possesses a distinguished basis
over $\QQ[a]$ formed by the
{\it double Schur functions\/} $s_{\la}(x\vt a)$ parameterized
by partitions $\la$.
The double Schur functions $s_{\la}(x\vt a)$
are closely related to the `factorial'
or `double' Schur polynomials $s_{\la}(x|a)$ which
were introduced by
Goulden and Greene~\cite{gg:nt} and Macdonald~\cite{m:sf}
as a generalization of the factorial Schur polynomials
of Biedenharn and Louck~\cite{bl:nc, bl:ib}.
Moreover, the polynomials $s_{\la}(x|a)$
are also obtained as a special case of
the double Schubert polynomials of Lascoux and
Sch\"utzenberger; see \cite{cll:fd}, \cite{l:sf}.
A formal definition of the ring $\La(x\vt a)$
and its basis elements $s_{\la}(x\vt a)$
can be found in a paper of Okounkov~\cite[Remark~2.11]{o:ni}
and reproduced below in Section~\ref{sec:def}.
The ring $\La$ is obtained from $\La(x\vt a)$ in
the specialization $a_i=0$ for all $i\in\ZZ$
while the elements $s_{\la}(x\vt a)$ turn into
the classical Schur functions $s_{\la}(x)\in\La$;
see Macdonald~\cite{m:sfh} for a detailed account
of the properties of $\La$.

Another specialization $a_i=-i+1$
for all $i\in\ZZ$ yields
the ring of {\it shifted symmetric functions\/}
$\La^*$, introduced
and studied by Okounkov and Olshanski~\cite{oo:ss}.
Many combinatorial results of \cite{oo:ss}
can be reproduced for the ring $\La(x\vt a)$ in a
rather straightforward way.
The respective specializations of the double Schur functions
in $\La^*$, known as the {\it shifted Schur functions\/}
were studied in \cite{o:qi}, \cite{oo:ss} in relation with
the higher Capelli identities
and quantum immanants for the Lie algebra $\gl_n$.

In a different kind of specialization,
the double Schur functions become
the equivariant Schubert classes
on Grassmannians; see e.g.
Knutson and Tao~\cite{kt:pe},
Fulton~\cite{f:ec} and
Mihalcea~\cite{m:gf}.
The structure coefficients $c_{\la\mu}^{\ts\nu}(a)$
of $\La(x\vt a)$ in the basis of
$s_{\la}(x\vt a)$, defined by the expansion
\beql{lrpoldef}
s_{\la}(x\vt a)\ts s_{\mu}(x\vt a)=
\sum_{\nu} c_{\la\mu}^{\ts\nu}(a)\ts s_{\tss\nu}(x\vt a),
\eeq
were called the {\it Littlewood--Richardson
polynomials\/} in \cite{m:lr}.
Under the respective specializations
they describe the multiplicative structure
of the equivariant cohomology ring
on the Grassmannian and the center
of the enveloping algebra $\U(\gl_n)$.
The polynomials $c_{\la\mu}^{\ts\nu}(a)$
possess the Graham positivity property:
they are polynomials in the differences $a_i-a_j$, $i<j$,
with positive integer coefficients; see \cite{g:pe}.
Explicit positive formulas for the polynomials $c_{\la\mu}^{\ts\nu}(a)$
were found in \cite{kt:pe}, \cite{k:elr} and \cite{m:lr};
an earlier formula found in \cite{ms:lr} lacks
the positivity property.
The Graham positivity brings natural combinatorics
of polynomials into the structure theory of $\La(x\vt a)$.
Namely, the entries of some transition
matrices between bases of $\La(x\vt a)$ such as
analogues of the Kostka numbers, turn out to be
Graham positive.

The {\it comultiplication\/} on the ring $\La(x\vt a)$ is
the $\QQ[a]$-linear ring homomorphism
\ben
\Delta:\La(x\vt a)\to \La(x\vt a)\ot^{}_{\ts\QQ[a]} \La(x\vt a)
\een
defined on the generators by
\ben
\Delta\big(p_k(x\vt a)\big)= p_k(x\vt a)\ot 1+1\ot p_k(x\vt a).
\een
In the specialization $a_i=0$ this homomorphism turns into the
comultiplication on the ring of symmetric functions $\La$;
see \cite[Chapter~I]{m:sfh}.
Define the {\it dual Littlewood--Richardson
polynomials\/}
$\wh c_{\la\mu}^{\ts\tss\nu}(a)$
as the coefficients in the expansion
\ben
\Delta\big(s_{\nu}(x\vt a)\big)=\sum_{\la,\ts\mu}
\wh c_{\la\mu}^{\ts\tss\nu}(a)\ts s_{\la}(x\vt a)
\ot s_{\mu}(x\vt a).
\een
The central problem we address in this paper is
calculation of the polynomials $\wh c_{\la\mu}^{\ts\tss\nu}(a)$
in an explicit form. Note that
if $|\nu|=|\la|+|\mu|$ then
$c_{\la\mu}^{\ts\nu}(a)=\wh c_{\la\mu}^{\ts\tss\nu}(a)
=c_{\la\mu}^{\ts\tss\nu}$ is the
Littlewood--Richardson coefficient.
Moreover,
\ben
c_{\la\mu}^{\ts\nu}(a)=0\quad\text{unless}
\quad |\nu|\leqslant|\la|+|\mu|,\fand
\wh c_{\la\mu}^{\ts\tss\nu}(a)=0\quad\text{unless}
\quad |\nu|\geqslant|\la|+|\mu|.
\een
We will show that the polynomials
$\wh c_{\la\mu}^{\ts\tss\nu}(a)$ can be interpreted
as the multiplication coefficients for certain
analogues of the Schur functions,
\ben
\wh s_{\la}(x\vt a)\ts \wh s_{\mu}(x\vt a)=
\sum_{\nu} \wh c_{\la\mu}^{\ts\tss\nu}(a)\ts \wh s_{\tss\nu}(x\vt a),
\een
where the $\wh s_{\la}(x\vt a)$ are symmetric
functions in $x$ which we call the {\it dual Schur functions\/}
(apparently, the term `dual double Schur functions'
would be more precise; we have chosen a shorter name
for the sake of brevity).
They can be given by the combinatorial formula
\beql{combi}
\wh s_{\la}(x\vt a)=\sum_{T}
\prod_{\al\in\la}X_{T(\al)}(a^{}_{-c(\al)+1},a^{}_{-c(\al)}),
\eeq
summed over the {\it reverse\/} $\la$-tableaux $T$, where
\ben
X_i(g,h)=
\frac{x_i\ts(1-g\ts x_{i-1})\dots (1-g\ts x_1)}
{(1-h\ts x_{i})\dots (1-h\ts x_1)},
\een
and $c(\al)=j-i$ denotes the content of the box $\al=(i,j)$;
see Section~\ref{sec:dua} below.

We calculate in an explicit form the
coefficients of the expansion of $\wh s_{\la}(x\vt a)$
as a series of the Schur functions $s_{\mu}(x)$ and vice versa.
This makes it possible to express
$\wh c_{\la\mu}^{\ts\tss\nu}(a)$ explicitly as polynomials
in the $a_i$ with the use of the Littlewood--Richardson
coefficients $c_{\la\mu}^{\ts\tss\nu}$.

The combinatorial formula \eqref{combi}
can be used to define the skew dual Schur functions, and
we show that the following decomposition holds
\ben
\wh s_{\tss\nu/\mu}(x\vt a)=
\sum_{\la} c_{\la\mu}^{\ts\nu}(a)\ts \wh s_{\la}(x\vt a),
\een
where the $c_{\la\mu}^{\ts\nu}(a)$ are the Littlewood--Richardson
polynomials.

The functions $\wh s_{\la}(x\vt a)$ turn out to be dual
to the double Schur functions
via the following analogue of the classical Cauchy identity:
\beql{cauchy}
\prod_{i,\ts j\geqslant 1}\frac{1-a_i\ts y_j}{1-x_i\ts y_j}
=\sum_{\la\in\Pc} s_{\la}(x\vt a)\ts \wh s_{\la}(y\vt a),
\eeq
where $\Pc$ denotes the set of all partitions and
$y=(y_1,y_2,\dots)$ is a set of variables.

The dual Schur functions $\wh s_{\la}(x\vt a)$
are elements of the extended
ring $\wh\La(x\vt a)$ of formal series of elements of
$\La(x)$ whose coefficients are polynomials in the $a_i$.
If $x=(x_1,x_2,\dots, x_n)$ is a finite set of variables
(i.e., $x_i=0$ for $i\geqslant n+1$), then $\wh s_{\la}(x\vt a)$
can be defined as the ratio
of alternants by analogy with the classical Schur polynomials.
With this definition of the dual Schur functions,
the identity \eqref{cauchy} can be deduced from
the `dual Cauchy formula' obtained in \cite[(6.17)]{m:sf}
and which is a
particular case of the Cauchy identity
for the double Schubert polynomials \cite{l:cc}.
An independent proof of a version of \eqref{cauchy}
for the shifted Schur functions (i.e., in the specialization
$a_i=-i+1$) was given by Olshanski~\cite{o:un}.
In the specialization $a_i=0$ each $\wh s_{\la}(x\vt a)$ becomes
the Schur function $s_{\la}(x)$, and \eqref{cauchy} turns into
the classical Cauchy identity.

We will also need a super version of the ring
of symmetric functions. The elements
\beql{poss}
p_k(x/y)=\sum_{i=1}^{\infty}\big(x_i^k+(-1)^{k-1}y_i^k\big)
\eeq
with $k=1,2,\dots$ are generators of
the ring of {\it supersymmetric functions\/}
which we will regard as a $\QQ[a]$-module
and denote by $\La(x/y\vt a)$. A distinguished basis
of $\La(x/y\vt a)$ was introduced by
Olshanski, Regev and Vershik~\cite{orv:fs}.
In a certain specialization the basis elements
become the {\it Frobenius--Schur functions\/}
$Fs_{\la}$ associated with the relative dimension
function on partitions; see \cite{orv:fs}.
In order to indicate dependence on the variables,
we will denote the basis elements by $s_{\la}(x/y\vt a)$
and call them the ({\it multiparameter\/})
{\it supersymmetric Schur functions\/}.
They
are closely related to the {\it factorial
supersymmetric Schur polynomials\/} introduced in \cite{m:fs};
see Section~\ref{sec:def} for precise formulas.
Note that the
evaluation map $y_i\mapsto -a_i$ for all $i\geqslant 1$
defines an isomorphism
\beql{isom}
\La(x/y\vt a)\to \La(x\vt a).
\eeq
The images of the generators \eqref{poss}
under this isomorphism
are the
double power sums symmetric functions
\eqref{powess}.
We will show that under the isomorphism \eqref{isom} we have
\beql{isomimsf}
s_{\la}(x/y\vt a)\mapsto s_{\la}(x\vt a).
\eeq
Due to \cite{orv:fs}, the
supersymmetric Schur functions possess a remarkable
combinatorial presentation in terms of diagonal-strict or
`shuffle' tableaux.
The isomorphism \eqref{isom}
implies the corresponding
combinatorial presentation for $s_{\la}(x\vt a)$
and allows us to introduce the {\it skew
double Schur functions\/}
$s_{\nu/\mu}(x\vt a)$.
The dual Littlewood--Richardson
polynomials
$\wh c_{\la\mu}^{\ts\tss\nu}(a)$ can then be
found from the expansion
\beql{skewde}
s_{\tss\nu/\mu}(x\vt a)=\sum_{\la}
\wh c_{\la\mu}^{\ts\tss\nu}(a)\ts s_{\la}(x\vt a),
\eeq
which leads to an alternative rule for the calculation of
$\wh c_{\la\mu}^{\ts\tss\nu}(a)$. This rule relies on
the combinatorial objects called `barred tableaux'
which were introduced in \cite{ms:lr}
for the calculation of the
polynomials $c_{\la\mu}^{\tss\nu}(a)$; see also
\cite{k:elr}, \cite{k:pf} and \cite{m:lr}.

The coefficients in the expansion of $s_{\mu}(x)$
in terms of the $\wh s_{\la}(x\vt a)$ turn out to coincide with
those in the decomposition of
$s_{\la}(x/y\vt a)$ in terms of the ordinary supersymmetric
Schur functions $s_{\la}(x/y)$ thus providing another
expression for these coefficients; cf. \cite{orv:fs}.

The identity \eqref{cauchy} allows us to introduce a pairing
between the rings $\La(x\vt a)$ and $\wh\La(x\vt a)$ so that
the respective families
$\{s_{\la}(x\vt a)\}$ and $\{\wh s_{\la}(x\vt a)\}$
are dual to each other. This leads to a natural definition
of the monomial and forgotten symmetric
functions in $\La(x\vt a)$ and $\wh\La(x\vt a)$
by analogy with \cite{m:sfh} and provides a relationship
between the transition matrices relating different
bases of these rings.

It is well known that
the ring of symmetric functions $\La$ admits an involutive
automorphism $\om:\La\to\La$ which interchanges the elementary
and complete symmetric functions; see \cite{m:sfh}.
We show that there is an isomorphism
$\om_a:\La(x\vt a)\to\La(x\vt a')$, and
$\om_a$ has the property $\om_{a'}\circ\om_{a}=\text{id}$,
where $a'$ denotes the sequence of parameters with $(a')_i=-a_{-i+1}$.
Moreover, the images of the natural bases elements of $\La(x\vt a)$
with respect to $\om_a$
can be explicitly described; see also \cite{oo:ss}
where such an involution was
constructed for the specialization $a_i=-i+1$, and
\cite{orv:fs} for its super version.
Furthermore, using a symmetry property of the
supersymmetric Schur functions,
we derive the symmetry properties of the
Littlewood--Richardson polynomials and their dual counterparts
\ben
c_{\la\mu}^{\tss\nu}(a)=c_{\la'\mu'}^{\tss\nu^{\tss\prime}}(a')
\Fand
\wh c_{\la\mu}^{\ts\tss\nu}(a)=
\wh c_{\la'\mu'}^{\ts\tss\nu^{\tss\prime}}(a'),
\een
where $\rho^{\tss\prime}$ denotes the conjugate
partition to any partition $\rho$.
In the context of equivariant cohomology,
the first relation is a consequence of the Grassmann duality;
see e.g. \cite[Lecture~8]{f:ec} and \cite{kt:pe}.

An essential
role in the proof of \eqref{cauchy} is played by
interpolation formulas for symmetric functions. The interpolation
approach goes back to the work of Okounkov~\cite{o:qi, o:ni},
where the key {\it vanishing theorem\/} for the double Schur functions
$s_{\la}(x\vt a)$ was proved; see also \cite{oo:ss}.
In a more general context, the Newton interpolation
for polynomials in several variables relies on the theory
of Schubert polynomials of Lascoux and
Sch\"utzenberger; see \cite{l:sf}.
The interpolation approach leads to a recurrence
relation for the coefficients $c_{P,\ts\mu}^{\ts\nu}(a)$
in the expansion
\beql{interp}
P\ts s_{\mu}(x\vt a)=\sum_{\nu}
c_{P,\ts\mu}^{\ts\nu}(a)\ts s_{\nu}(x\vt a),\qquad P\in\La(x\vt a),
\eeq
as well as to an explicit formula for the $c_{P,\ts\mu}^{\ts\nu}(a)$
in terms of the values of $P$; see \cite{ms:lr}.
Therefore, the (dual) Littlewood--Richardson polynomials
and the entries of the transition matrices between
various bases of $\La(x\vt a)$ can be given as rational functions
in the variables $a_i$.
Under appropriate specializations,
these formulas imply some
combinatorial identities
involving Kostka numbers, irreducible
characters of the symmetric
group and dimensions of skew diagrams; cf. \cite{oo:ss}.

\medskip

I am grateful to Grigori Olshanski for valuable
remarks and discussions.

\section{Double and supersymmetric Schur functions}
\label{sec:def}
\setcounter{equation}{0}

\subsection{Definitions and preliminaries}
\label{subsec:dpr}

Recall the definition of the ring $\La(x\vt a)$ from
\cite[Remark~2.11]{o:ni};
see also \cite{m:lr}.
For each nonnegative integer $n$
denote by $\La_n$ the ring of symmetric
polynomials in $x_1,\dots,x_n$ with coefficients in $\QQ[a]$
and let $\La^k_n$ denote the $\QQ[a]$-submodule of $\La_n$
which consists of
the polynomials $P_n(x_1,\dots,x_n)$ such that the total degree
of $P_n$ in the variables $x_i$ does not exceed $k$.
Consider the evaluation maps
\beql{eval}
\varphi_n:\La^{k}_{n}\to\La^{k}_{n-1},\qquad
P_n(x_1,\dots,x_n)\mapsto
P_n(x_1,\dots,x_{n-1},a_n)
\eeq
and the corresponding inverse limit
\ben
\La^{k}=\lim_{\longleftarrow}
\La^{k}_n,\qquad n\to\infty.
\een
The elements of $\La^{k}$ are sequences
$P=(P_0,P_1,P_2,\dots)$ with $P_n\in \La^{k}_n$
such that
\ben
\varphi_n(P_n)=P_{n-1}
\qquad\text{for}\quad n=1,2,\dots.
\een
Then the union
\ben
\La(x\vt a)= \bigcup_{k\geqslant 0}\La^{k}
\een
is a ring with the product
\ben
P\tss Q=(P_0\tss Q_0,P_1\tss Q_1,P_2\tss Q_2,\dots),
\qquad Q=(Q_0,Q_1,Q_2,\dots).
\een
The elements of $\La(x\vt a)$ may be regarded
as formal series in the variables
$x_i$ with coefficients in $\QQ[a]$.
For instance, the sequence of polynomials
\ben
\sum_{i=1}^n (x_i^k-a_i^k),\qquad n \geqslant 0,
\een
determines the {\it double power sums symmetric function\/}
\eqref{powess}.

Note that if $k$ is fixed, then the evaluation
maps \eqref{eval} are isomorphisms for
all sufficiently large values of $n$. This allows one
to establish many properties of $\La(x\vt a)$
by working with finite sets of variables $x=(x_1,\dots,x_n)$.

Now we recall the definition and some key properties of the
double Schur functions. We basically follow
\cite[6th~Variation]{m:sf} and \cite{o:ni},
although our notation is slightly different.
A partition
$\la$ is a weakly decreasing sequence
$\la=(\la_1,\dots,\la_l)$ of integers $\la_i$ such that
$\la_1\geqslant\dots\geqslant\la_l\geqslant 0$.
Sometimes this sequence
is considered to be completed by a finite or infinite
sequence of zeros.
We will identify $\la$ with its diagram
represented graphically as
the array of left justified rows of unit boxes
with $\la_1$ boxes in the top row, $\la_2$ boxes in the second
row, etc. The total number of boxes in $\la$
will be denoted by $|\la|$ and the number of nonzero rows
will be called the length of $\la$ and denoted
$\ell(\la)$.
The transposed diagram
$\la'=(\la'_1,\dots,\la'_p)$ is obtained from $\la$
by applying the symmetry
with respect to the main diagonal, so that
$\la'_j$ is the number of boxes in the $j$-th column of $\la$.
If $\mu$ is a diagram contained in $\la$, then the
skew diagram $\la/\mu$ is the set-theoretical difference
of diagrams $\la$ and $\mu$.

Suppose now that $x=(x_1,\dots,x_n)$ is a finite set
of variables. For any $n$-tuple of nonnegative integers
$\al=(\al_1,\dots,\al_n)$ set
\ben
A_{\al}(x\vt a)=\det\big[(x_i\vt a)^{\al_j}\big]_{i,j=1}^n,
\een
where $(x_i\vt a)^0=1$ and
\ben
(x_i\vt a)^r=(x_i-a_n)(x_i-a_{n-1})\dots (x_i-a_{n-r+1}),
\qquad r\geqslant 1.
\een
For any partition $\la=(\la_1,\dots,\la_n)$ of length
not exceeding $n$
set
\ben
s_{\la}(x\vt a)=\frac{A_{\la+\de}(x\vt a)}{A_{\de}(x\vt a)},
\een
where $\de=(n-1,\dots,1,0)$. Note that since $A_{\de}(x\vt a)$
is a skew-symmetric polynomial in $x$ of degree $n(n-1)/2$, it
coincides with the Vandermonde determinant,
\ben
A_{\de}(x\vt a)=\prod_{1\leqslant i<j\leqslant n}(x_i-x_j)
\een
and so $s_{\la}(x\vt a)$ belongs to the ring $\La_n$.
Moreover,
\ben
s_{\la}(x\vt a)=s_{\la}(x)+\ \text{lower degree terms in}\ \  x,
\een
where $s_{\la}(x)$ is the Schur polynomial; see e.g.
\cite[Chapter~I]{m:sfh}.
We also set $s_{\la}(x\vt a)=0$ if $\ell(\la)> n$.
Then under the evaluation map \eqref{eval} we have
\ben
\varphi_n: s_{\la}(x\vt a)\mapsto s_{\la}(x'\vt a),
\qquad x'=(x_1,\dots,x_{n-1}),
\een
so that the sequence $\big(s_{\la}(x\vt a)\ |\ n\geqslant 0\big)$
defines an element of the ring $\La(x\vt a)$. We will keep
the notation $s_{\la}(x\vt a)$ for this element of
$\La(x\vt a)$, where $x$ is now understood as the infinite sequence of
variables, and call it the {\it double Schur function\/}.

By a {\it reverse $\la$-tableau\/}
$T$ we will mean a tableau obtained by filling in the boxes
of $\la$ with the positive integers in such a way that
the entries weakly decrease along the rows and strictly decrease
down the columns. If $\al=(i,j)$ is a box of $\la$
in row $i$ and column $j$,
we let $T(\al)=T(i,j)$ denote the entry of $T$ in the box $\al$
and let $c(\alpha)=j-i$ denote the content of
this box.
The double Schur functions admit the following
tableau presentation
\beql{defdouble}
s_{\la}(x\vt a)=\sum_{T}
\prod_{\al\in\la}
(x^{}_{T(\alpha)}-a^{}_{T(\alpha)-c(\al)}),
\eeq
summed over all reverse $\la$-tableaux $T$.

When the entries of $T$ are restricted to
the set $\{1,\dots,n\}$, formula \eqref{defdouble}
provides the respective tableau presentation of
the polynomials $s_{\la}(x\vt a)$ with $x=(x_1,\dots,x_n)$.
Moreover, in this case the formula can be extended to
skew diagrams and we define the corresponding
polynomials by
\beql{defdoubleskew}
\wt s_{\theta}(x\vt a)=\sum_{T}
\prod_{\al\in\theta}
(x^{}_{T(\alpha)}-a^{}_{T(\alpha)-c(\al)}),
\eeq
summed over all reverse $\theta$-tableaux $T$
with entries in $\{1,\dots,n\}$, where $\theta$
is a skew diagram. We suppose that $\wt s_{\theta}(x\vt a)=0$
unless all columns of $\theta$ contain at most $n$ boxes.

\bre\label{rem:skew}
(i)\quad
Although the polynomials \eqref{defdoubleskew}
belong to the ring $\La_n$, they are generally
not consistent with respect to
the evaluation maps \eqref{eval}.
We used different notation in \eqref{defdouble} and
\eqref{defdoubleskew} in order to distinguish between
the polynomials $\wt s_{\theta}(x\vt a)$ and the
skew double Schur functions $s_{\theta}(x\vt a)$
to be introduced in
Definition~\ref{def:skdsf} below.

\medskip
\noindent
(ii)\quad
In order to relate our notation to \cite{m:sf}, note that
for the polynomials $\wt s_{\theta}(x\vt a)$ with $x=(x_1,\dots,x_n)$
we have
\ben
\wt s_{\theta}(x\vt a)=s_{\theta}(x\tss|\tss u),
\een
where the sequences $a=(a_i)$ and $u=(u_i)$ are related by
\beql{seqs}
u_i=a_{n-i+1},\qquad i\in\ZZ.
\eeq
The polynomials $s_{\theta}(x\tss|\tss u)$ are often called
the {\it factorial Schur polynomials\/} ({\it functions\/})
in the literature. They can be given by
the combinatorial formula
\beql{facsf}
s_{\theta}(x\tss|\tss u)=\sum_{T}
\prod_{\al\in\theta}
(x^{}_{T(\alpha)}-u^{}_{T(\alpha)+c(\al)}),
\eeq
summed over all {\it semistandard\/} $\theta$-tableaux $T$
with entries in $\{1,\dots,n\}$; the entries of $T$
weakly increase along the rows and strictly increase
down the columns.

\medskip
\noindent
(iii)\quad
If we replace $a_i$ with $c_{-i}$
and index the variables $x$ with nonnegative integers,
the double Schur functions $s_{\la}(x\vt a)$ will become
the corresponding symmetric
functions of \cite{o:ni}; cf. formula (3.7) in that paper.
Moreover, under the specialization $a_i=-i+1$ for all $i\in\ZZ$
the double Schur functions become the {\it shifted
Schur functions\/} of \cite{oo:ss} in the variables $y_i=x_i+i-1$.
\qed
\ere

\subsection{Analogues of classical bases}

The {\it double elementary\/} and {\it complete
symmetric functions} are defined respectively by
\ben
e_k(x\vt a)=s_{(1^k)}(x\vt a),\qquad h_k(x\vt a)=s_{(k)}(x\vt a)
\een
and hence, they can be given by the formulas
\ben
\bal
e_k(x\vt a)&=\sum_{i_1>\dots>i_k}(x_{i_1}-a_{i_1})\dots
(x_{i_k}-a_{i_k+k-1}),\\
h_k(x\vt a)&=\sum_{i_1\geqslant\dots\geqslant i_k}
(x_{i_1}-a_{i_1})\dots
(x_{i_k}-a_{i_k-k+1}).
\eal
\een
Their generating functions can be written
by analogy with the classical case as in \cite{m:sfh}
and they take the form
\begin{align}\label{gene}
1+\sum_{k=1}^{\infty} \frac{e_k(x\vt a)\ts t^k}
{(1+a_1\tss t)\dots(1+a_k\tss t)}&=
\prod_{i=1}^{\infty}\frac{1+x_i\tss t}{1+a_i\tss t},\\[1em]
\label{genh}
1+\sum_{k=1}^{\infty} \frac{h_k(x\vt a)\ts t^k}
{(1-a_0\tss t)\dots(1-a_{-k+1}\tss t)}&=
\prod_{i=1}^{\infty}\frac{1-a_i\tss t}{1-x_i\tss t};
\end{align}
see e.g. \cite{m:sf}, \cite{oo:ss}.

Given a partition
$\la=(\la_1,\dots,\la_l)$, set
\ben
\bal
p_{\la}(x\vt a)&=p_{\la_1}(x\vt a)\dots p_{\la_l}(x\vt a),\\
e_{\la}(x\vt a)&=e_{\la_1}(x\vt a)\dots e_{\la_l}(x\vt a),\\
h_{\la}(x\vt a)&=h_{\la_1}(x\vt a)\dots h_{\la_l}(x\vt a).
\eal
\een

The following proposition is easy to deduce from
the properties of the classical symmetric functions;
see \cite{m:sfh}.

\bpr\label{prop:basis}
Each of the families $p_{\la}(x\vt a)$,
$e_{\la}(x\vt a)$, $h_{\la}(x\vt a)$
and $s_{\la}(x\vt a)$, parameterized by all
partitions $\la$,
forms a basis of $\La(x\vt a)$ over $\QQ[a]$.
\qed
\epr

In particular, each of the families
$p_k(x\vt a)$,
$e_k(x\vt a)$ and $h_k(x\vt a)$ with $k\geqslant 1$
is a set of algebraically independent generators of
$\La(x\vt a)$ over $\QQ[a]$.
Under the specialization $a_i=0$, the bases
of Proposition~\ref{prop:basis} turn into
the classical bases $p_{\la}(x)$,
$e_{\la}(x)$, $h_{\la}(x)$
and $s_{\la}(x)$ of $\La$.
The ring of symmetric functions $\La$
possesses two more bases $m_{\la}(x)$ and $f_{\la}(x)$;
see \cite[Chapter~I]{m:sfh}.
The {\it monomial
symmetric functions\/} $m_{\la}(x)$ are defined by
\ben
m_{\la}(x)=\sum_{\si} x_{\si(1)}^{\la_1}x_{\si(2)}^{\la_2}\dots
x_{\si(l)}^{\la_l},
\een
summed over permutations $\si$ of the $x_i$
which give distinct monomials. The basis elements $f_{\la}(x)$
are called the {\it forgotten symmetric functions\/}, they are
defined as the images of the $m_{\la}(x)$ under the involution
$\om:\La\to\La$ which takes $e_{\la}(x)$ to $h_{\la}(x)$;
see \cite{m:sfh}.
The corresponding basis elements
$m_{\la}(x\vt a)$ and $f_{\la}(x\vt a)$
in $\La(x\vt a)$ will be defined in Section~\ref{sec:ome}.

\subsection{Duality isomorphism}

Introduce the sequence of variables $a'$ which is
related to the sequence $a$ by the rule
\ben
(a')_i=-a_{-i+1},\qquad i\in\ZZ.
\een
The operation $a\mapsto a'$ is clearly involutive so that
$(a')'=a$. Note that any element of
the polynomial ring $\QQ[a']$ can be identified with
the element of $\QQ[a]$ obtained by replacing each $(a')_i$
by $-a_{-i+1}$.
Define the ring homomorphism
\ben
\om_a:\La(x\vt a)\to\La(x\vt a')
\een
as the $\QQ[a]$-linear map
such that
\beql{omega}
\om_a: e_k(x\vt a)\mapsto h_k(x\vt a'),\qquad k=1,2,\dots.
\eeq
An arbitrary element of $\La(x\vt a)$ can be written
as a unique linear combination of the basis elements $e_{\la}(x\vt a)$
with coefficients in $\QQ[a]$. The image of such
a linear combination under $\om_a$ is then found by
\ben
\om_a:\sum_{\la}c_{\la}(a)\ts e_{\la}(x\vt a)
\mapsto \sum_{\la}c_{\la}(a)\ts h_{\la}(x\vt a'),\qquad
c_{\la}(a)\in\QQ[a],
\een
and $c_{\la}(a)$ is regarded as an element of $\QQ[a']$.
Clearly, $\om_a$ is a ring isomorphism,
since the $h_k(x\vt a')$ are
algebraically independent generators of
$\La(x\vt a')$ over $\QQ[a']$. In the case of finite set
of variables $x=(x_1,\dots,x_n)$ the respective
isomorphism $\om_a$ is defined by the same rule \eqref{omega}
with the values $k=1,\dots,n$.

\bpr\label{prop:om}
We have
$\om_{a'}\circ \om_a={\rm id}^{}_{\La(x\vt a)}$ and
\beql{omhla}
\om_a: h_{\la}(x\vt a)\mapsto e_{\la}(x\vt a').
\eeq
\epr

\bpf
Relations \eqref{gene} and \eqref{genh} imply that
\ben
\Bigg(\sum_{k=0}^{\infty} \frac{(-1)^k\ts e_k(x\vt a)\ts t^k}
{(1-a_1\tss t)\dots(1-a_k\tss t)}\Bigg)
\Bigg(\sum_{r=0}^{\infty} \frac{h_r(x\vt a)\ts t^r}
{(1-a_0\tss t)\dots(1-a_{-r+1}\tss t)}\Bigg)=1.
\een
Applying the isomorphism $\om_a$, we get
\ben
\Bigg(\sum_{k=0}^{\infty} \frac{(-1)^k\ts h_k(x\vt a')\ts t^k}
{(1+(a')_0\tss t)\dots(1+(a')_{-k+1}\tss t)}\Bigg)
\Bigg(\sum_{r=0}^{\infty} \frac{\om_a\big(h_r(x\vt a)\big)\ts t^r}
{(1+(a')_1\tss t)\dots(1+(a')_r\tss t)}\Bigg)=1.
\een
Replacing here $t$ by $-t$ and comparing with the previous
identity, we can conclude that
$\om_a\big(h_r(x\vt a)\big)=e_r(x\vt a')$.
This proves \eqref{omhla} and the first part
of the proposition, because
$\om_{a'}\big(h_r(x\vt a')\big)=e_r(x\vt a)$.
\epf

We will often use the shift operator
$\tau$ whose
powers act on sequences by the rule
\ben
(\tau^k a)_i=a_{k+i}\qquad\text{for}\quad k\in\ZZ.
\een
The following analogues of the
Jacobi--Trudi
and N\"{a}gelsbach--Kostka formulas
are immediate from \cite[(6.7)]{m:sf}.
Namely, if the set of variables $x=(x_1,\dots,x_n)$ is
finite and $\la$
is a partition of length not exceeding $n$, then
\beql{jt}
s_{\la}(x\vt a)
=\det\big[h_{\la_i-i+j}(x\vt \tau^{\tss j-1}\tss a)\big]
\eeq
and
\beql{nk}
s_{\la}(x\vt a)
=\det\big[e_{\la'_i-i+j}(x\vt \tau^{\tss -j+1}\tss a)\big],
\eeq
where the determinants are taken over the respective sets of indices
$i,j=1,\dots,\ell(\la)$ and
$i,j=1,\dots,\ell(\la')$.

\subsection{Skew double Schur functions}

Consider now the ring of supersymmetric functions
$\La(x/y\vt a)$ defined in the Introduction.
Taking two finite sets of variables $x=(x_1,\dots,x_n)$
and $y=(y_1,\dots,y_n)$, define the {\it supersymmetric
Schur polynomial\/} $s_{\nu/\mu}(x/y\vt a)$
associated with a skew diagram $\nu/\mu$ by the formula
\beql{defssf}
s_{\nu/\mu}(x/y\vt a)=\sum_{\mu\subseteq\ts\rho\ts\subseteq\tss\nu}
\wt s_{\nu/\rho}(x\vt a)\ts
s_{\rho^{\tss\prime}/\mu^{\tss\prime}}(y \tss|\tss {-}a),
\eeq
where the polynomials $\wt s_{\nu/\rho}(x\vt a)$ and
$s_{\rho^{\tss\prime}/\mu^{\tss\prime}}(y \tss|\tss{-}a)$ are defined
by the respective combinatorial formulas
\eqref{defdoubleskew} and \eqref{facsf}.
The polynomials \eqref{defssf} coincide with
the factorial supersymmetric Schur polynomials
$s_{\nu/\mu}(x/y\tss|\tss u)$
of \cite{m:fs} associated with the sequence
$u$ related to $a$ by \eqref{seqs}.
It was observed in \cite{orv:fs}
that the sequence of polynomials
$\big(s_{\nu/\mu}(x/y\vt a)\ts|\ts n\geqslant 1\big)$ is consistent
with respect to the evaluations $x_n=y_n=0$ and hence,
it defines the {\it supersymmetric
Schur function\/} $s_{\nu/\mu}(x/y\vt a)$, where
$x$ and $y$ are infinite sequences of variables
(in fact, Proposition~3.4 in \cite{orv:fs} needs to be
extended to skew diagrams which is immediate).
Moreover, in \cite{orv:fs}
these functions were given by new
combinatorial formulas.
In order to write them down,
consider the ordered alphabet
\ben
\AAb=\{1'<1<2^{\tss\prime}<2<\dots\}.
\een
Given a skew diagram $\theta$, an
$\AAb$-{\it tableau\/} $T$ of shape $\theta$
is obtained by filling in the boxes of $\theta$ with
the elements of $\AAb$ in such a way that the
entries of $T$ weakly increase along each row and
down each column, and for each $i=1,2,\dots$ there is at most
one symbol $i'$ in each row and at most
one symbol $i$ in each column of $T$.
The following formula gives
the supersymmetric
Schur function $s_{\theta}(x/y\vt a)$ associated with $\theta$:
\beql{tabapr}
s_{\theta}(x/y\vt a)=
\sum_{T}
\prod_{\underset{\scriptstyle T(\alpha)
\text{\ts\ts unprimed}}{\alpha\in\theta}}
\big(x^{}_{T(\al)}-a_{-c(\al)+1}\big)
\prod_{\underset{\scriptstyle T(\alpha)
\text{\ts\ts primed}}{\alpha\in\theta}}
\big(y^{}_{T(\al)}+a_{-c(\al)+1}\big),
\eeq
summed over all $\AAb$-tableaux $T$ of shape $\theta$,
where the subscripts of the variables $y_i$ are identified
with the primed indices. An alternative formula is obtained
by using a different ordering of the alphabet:
\ben
\AAb'=\{1<1'<2<2^{\tss\prime}<\dots\}.
\een
The $\AAb'$-tableaux $T$ of shape $\theta$ are defined
in exactly the same way as the $\AAb$-tableaux, only taking
into account the new ordering. Then
\beql{taba}
s_{\theta}(x/y\vt a)=
\sum_{T}
\prod_{\underset{\scriptstyle T(\alpha)
\text{\ts\ts unprimed}}{\alpha\in\theta}}
\big(x^{}_{T(\al)}-a_{-c(\al)}\big)
\prod_{\underset{\scriptstyle T(\alpha)
\text{\ts\ts primed}}{\alpha\in\theta}}
\big(y^{}_{T(\al)}+a_{-c(\al)}\big),
\eeq
summed over all $\AAb'$-tableaux $T$ of shape $\theta$.

The supersymmetric
Schur functions have the following symmetry property
\beql{symprsu}
s_{\theta}(x/y\vt a)=s_{\theta^{\tss\prime}}(y/x\vt a')
\eeq
implied by their combinatorial presentation.
Moreover, if $x_i=y_i=0$ for all $i\geqslant n+1$, then
only tableaux $T$ with entries in $\{1,1',\dots,n,n'\}$
make nonzero contributions in
either \eqref{tabapr} or \eqref{taba}.

\bre\label{rem:deforv}
The supersymmetric
Schur function
$s_{\theta}(x/y\vt a)$ given in \eqref{tabapr}
coincides with $\Sigma_{\theta;-a'}(x;y)$ as defined
in \cite[Proposition~4.4]{orv:fs}.
In order to derive \eqref{taba}, first use \eqref{symprsu},
then apply the transposition of the tableaux
with respect to the main diagonal and swap $i$ and $i^{\tss\prime}$
for each $i$.
Note that \cite{orv:fs} also contains an
equivalent combinatorial formula for $\Sigma_{\theta;a}(x;y)$
in terms of skew hooks.
\qed
\ere

\bpr\label{prop:imisom}
The image of the supersymmetric
Schur function $s_{\nu}(x/y\vt a)$ associated with
a {\rm(}nonskew{\rm)} diagram $\nu$ under the isomorphism \eqref{isom}
coincides with the double Schur function
$s_{\nu}(x\vt a)$; that is,
\ben
s_{\nu}(x/y\vt a)\big|^{}_{y=-a}=s_{\nu}(x\vt a),
\een
where $y=-a$ denotes the evaluation $y_i=-a_i$ for $i\geqslant 1$.
\epr

\bpf
We may assume that the sets of variables $x$ and $y$
are finite,
$x=(x_1,\dots,x_n)$
and $y=(y_1,\dots,y_n)$. The claim
now follows from relation \eqref{defssf} with $\mu=\varnothing$,
if we observe that
$s_{\rho^{\tss\prime}}(y \tss|\tss {-}a)\big|^{}_{y=-a}=0$
unless $\rho=\varnothing$.
\epf

The symmetry property \eqref{symprsu} implies
the following dual version of Proposition~\ref{prop:imisom}.

\bco\label{cor:dueva}
Under the isomorphism $\La(x/y\vt a)\to \La(y\vt a')$
defined by the evaluation $x_i=-(a')_i$ for all $i\geqslant 1$
we have
\ben
s_{\theta}(x/y\vt a)\big|^{}_{x=-a'}=s_{\theta'}(y\vt a').
\een
\eco

Using Proposition~\ref{prop:imisom}, we can find the images
of the double Schur functions with respect to
the duality isomorphism $\om_a$ defined in \eqref{omega}.

\bco\label{cor:imsuf}
Under the isomorphism $\om_a:\La(x\vt a)\to\La(x\vt a')$
we have
\beql{omsla}
\om_a: s_{\la}(x\vt a)\mapsto s_{\la'}(x\vt a').
\eeq
\eco

\bpf
The Littlewood--Richardson
polynomials $c_{\la\mu}^{\tss\nu}(a)$ are defined
by the expansion \eqref{lrpoldef}. Hence,
by Proposition~\ref{prop:imisom}
we have
\ben
s_{\la}(x/y\vt a)\ts s_{\mu}(x/y\vt a)=
\sum_{\nu} c_{\la\mu}^{\ts\nu}(a)\ts s_{\tss\nu}(x/y\vt a).
\een
Using \eqref{symprsu}, we get
\beql{symlr}
c_{\la\mu}^{\tss\nu}(a)=c_{\la'\mu'}^{\tss\nu^{\tss\prime}}(a').
\eeq
Now, observe that relation \eqref{omsla} can be taken as
a definition of the $\QQ[a]$-module isomorphism
$\La(x\vt a)\to\La(x\vt a')$. Moreover, this definition
agrees with \eqref{omega}.
Therefore, it is sufficient to verify that this
$\QQ[a]$-module isomorphism is a ring homomorphism.
Applying \eqref{symlr} we obtain
\ben
\bal
\om_a\big(s_{\la}(x\vt a)\ts s_{\mu}(x\vt a)\big)
{}&{}=\sum_{\nu} c_{\la\mu}^{\ts\nu}(a)\ts
\om_a\big(s_{\tss\nu}(x\vt a)\big)
=\sum_{\nu} c_{\la'\mu'}^{\ts\nu^{\tss\prime}}(a')\ts
s_{\tss\nu^{\tss\prime}}(x\vt a')\\
{}&{}=s_{\la'}(x\vt a')\ts s_{\mu'}(x\vt a')=
\om_a\big(s_{\la}(x\vt a)\big)\ts \om_a\big(s_{\mu}(x\vt a)\big).
\eal
\een
\epf

Proposition~\ref{prop:imisom} leads to the following definition.

\bde\label{def:skdsf}
For any skew diagram $\theta$ define the
{\it skew double Schur function\/} $s_{\theta}(x\vt a)\in\La(x\vt a)$
as the image of $s_{\theta}(x/y\vt a)\in\La(x/y\vt a)$
under the isomorphism \eqref{isom}; that is,
\ben
s_{\theta}(x\vt a)=s_{\theta}(x/y\vt a)\big|^{}_{y=-a}.
\een
Equivalently, using \eqref{tabapr} and \eqref{taba},
respectively, we have
\beql{tabadsdu}
s_{\theta}(x\vt a)=
\sum_{T}
\prod_{\underset{\scriptstyle T(\alpha)
\text{\ts\ts unprimed}}{\alpha\in\theta}}
\big(x^{}_{T(\al)}-a_{-c(\al)+1}\big)
\prod_{\underset{\scriptstyle T(\alpha)
\text{\ts\ts primed}}{\alpha\in\theta}}
\big(a_{-c(\al)+1}-a^{}_{T(\al)}\big),
\eeq
summed over all $\AAb$-tableaux $T$ of shape $\theta$; and
\beql{tabads}
s_{\theta}(x\vt a)=
\sum_{T}
\prod_{\underset{\scriptstyle T(\alpha)
\text{\ts\ts unprimed}}{\alpha\in\theta}}
\big(x^{}_{T(\al)}-a_{-c(\al)}\big)
\prod_{\underset{\scriptstyle T(\alpha)
\text{\ts\ts primed}}{\alpha\in\theta}}
\big(a_{-c(\al)}-a^{}_{T(\al)}\big),
\eeq
summed over all $\AAb'$-tableaux $T$ of shape $\theta$.
Furthermore, by \eqref{defssf} the
skew double Schur function $s_{\nu/\mu}(x\vt a)$
can also be defined as the sequence of polynomials
\beql{seqpol}
s_{\nu/\mu}(x\vt a)=\sum_{\mu\subseteq\ts\rho\ts\subseteq\tss\nu}
\wt s_{\nu/\rho}(x\vt a)\ts
s_{\rho^{\tss\prime}/\mu^{\tss\prime}}(-a^{(n)} \tss|\tss {-}a),
\qquad n=1,2,\dots,
\eeq
where $x=(x_1,\dots,x_n)$ and $a^{(n)}=(a_1,\dots,a_n)$.
\qed
\ede

For any partition $\mu$ introduce the sequence $a_{\mu}$
and the series $|a_{\mu}|$ by
\ben
a_{\mu}=(a_{1-\mu_1},a_{2-\mu_2},\dots)\Fand
|a_{\mu}|=a_{1-\mu_1}+a_{2-\mu_2}+\dots.
\een
Given any element $P(x)\in\La(x\vt a)$, the value
$P(a_{\mu})$
is a well-defined element of $\QQ[a]$.
The vanishing theorem of Okounkov~\cite{o:qi, o:ni}
states that
\ben
s_{\la}(a_{\rho}\vt a)=0\qquad\text{unless}
\quad \la\subseteq\rho,
\een
and
\beql{hoo}
s_{\la}(a_{\la}\vt a)=\prod_{(i,j)\in\la}
\big(a^{}_{i-\la_i}-a^{}_{\la'_j-j+1}\big).
\eeq
This theorem can be used to derive the
interpolation formulas given in the next
proposition. In a slightly different situation this derivation
was performed in \cite[Propositions~3.3 \& 3.4]{ms:lr}
relying on the approach of \cite{oo:ss},
and an obvious
modification of those arguments works in the present context;
see also \cite{f:ec}, \cite{kt:pe}.
The expressions like $|a_{\nu}|-|a_{\mu}|$ used
below are understood
as the polynomials $\sum_{i\geqslant 1}
(a_{i-\nu_i}-a_{i-\mu_i})$.
We will write $\rho\to\sigma$ if the diagram $\si$
is obtained from the diagram $\rho$
by adding one box.

\bpr\label{prop:interp}
Given an element $P(x)\in\La(x\vt a)$, define the
polynomials $c_{P,\ts\mu}^{\ts\nu}(a)$
by the expansion
\beql{interpo}
P(x)\ts s_{\mu}(x\vt a)=\sum_{\nu}
c_{P,\ts\mu}^{\ts\nu}(a)\ts s_{\nu}(x\vt a).
\eeq
Then
$c_{P,\ts\mu}^{\ts\nu}(a)=0$ unless $\mu\subseteq \nu$,
and $c_{P,\ts\mu}^{\ts\mu}(a)=P(a_{\mu})$. Moreover, if
$\mu\subseteq \nu$, then
\ben
c_{P,\ts\mu}^{\ts\nu}(a)=\frac{1}{|a_{\nu}|-|a_{\mu}|}
\Bigg(\sum_{\mu^+,\ts\mu\to\mu^+}c_{P,\ts\mu^+}^{\ts\nu}(a)
-\sum_{\nu^-,\ts\nu^-\to\nu}c_{P,\ts\mu}^{\ts\nu^-}(a)\Bigg).
\een
The same coefficient can also be found by the
formula
\beql{rati}
c_{P,\ts\mu}^{\ts\nu}(a)=\sum_{R}\sum_{k=0}^{l}
\frac{P(a_{\rho^{(k)}})}{(|a_{\rho^{(k)}}|-|a_{\rho^{(0)}}|)
\ldots\wedge\ldots(|a_{\rho^{(k)}}|-|a_{\rho^{(l)}}|)},
\eeq
summed over all sequences of partitions $R$
of the form
\ben
\mu=\rho^{(0)}\to\rho^{(1)}\to
\dots\to\rho^{(l-1)}\to\rho^{(l)}=\nu,
\een
where the symbol $\wedge$ indicates that the zero factor
should be skipped.
\qed
\epr

\section{Cauchy identities and dual Schur functions}
\label{sec:dua}
\setcounter{equation}{0}

\subsection{Definition of dual Schur functions and
Cauchy identities}

We let $\wh\La(x\vt a)$ denote the ring
of formal series of the symmetric functions in
the set of indeterminates $x=(x_1,x_2,\dots)$
with coefficients in $\QQ[a]$. More precisely,
\beql{defla}
\wh\La(x\vt a)=\Big\{\sum_{\la\in\Pc} c_{\la}(a)\ts s_{\la}(x)
\ |\ c_{\la}(a)\in \QQ[a]\Big\}.
\eeq
The Schur functions $s_{\la}(x)$ can certainly
be replaced here by any other classical basis of $\La$ parameterized
by the set of partitions $\Pc$. We will use the symbol
$\wh\La_n=\wh\La_n(x\vt a)$ to indicate
the ring defined as in \eqref{defla}
for the case of the finite set of variables $x=(x_1,\dots,x_n)$.
An element of $\wh\La(x\vt a)$ can be viewed as a sequence
of elements of $\wh\La_n$ with $n=0,1,\dots$,
consistent with respect to the evaluation maps
\ben
\psi_n: \wh\La_n\to \wh\La_{n-1},\qquad
Q(x_1,\dots,x_n)\mapsto Q(x_1,\dots,x_{n-1},0).
\een

For any $n$-tuple of nonnegative integers
$\be=(\be_1,\dots,\be_n)$ set
\ben
A_{\be}(x, a)=\det\big[(x_i,a)^{\be_j}
\tss(1-a^{}_{n-\be_j-1}x_i)(1-a^{}_{n-\be_j-2}\ts x_i)\dots
(1-a^{}_{1-\be_j}x_i)\big]_{i,j=1}^n,
\een
where $(x_i,a)^{0}=1$ and
\beql{dumon}
(x_i,a)^{r}=\frac{x_i^r}
{(1-a^{}_0\tss x_i)(1-a^{}_{-1}\tss x_i)\dots (1-a^{}_{1-r}\tss x_i)},
\qquad r\geqslant 1.
\eeq

Let $\la=(\la_1,\dots,\la_n)$ be a partition of length
not exceeding $n$. Denote by $d$ the number of boxes
on the diagonal of $\la$. That is, $d$ is
determined by
the condition that
$\la_{d+1}\leqslant d\leqslant \la_d$.
The $(i,j)$ entry $A_{ij}$ of the determinant
$A_{\la+\de}(x, a)$ can be written
more explicitly as
\ben
A_{ij}=\begin{cases}\dfrac{x_i^{\la_j+n-j}}
{(1-a^{}_0\tss x_i)(1-a^{}_{-1}\tss x_i)\dots (1-a^{}_{j-\la_j}\tss x_i)}
\quad&\text{for}\quad j=1,\dots,d,\\[1.5em]
x_i^{\la_j+n-j}\ts
(1-a^{}_1\tss x_i)(1-a^{}_2\tss x_i)\dots (1-a^{}_{j-\la_j-1}\tss x_i)
\quad&\text{for}\quad j=d+1,\dots,n.
\end{cases}
\een
Observe that the determinant $A_{\de}(x, a)$
corresponding to the empty partition equals the Vandermonde
determinant,
\ben
A_{\de}(x, a)=\prod_{1\leqslant i<j\leqslant n}(x_i-x_j).
\een
Hence, the formula
\beql{dusf}
\wh s_{\la}(x\vt a)=\frac{A_{\la+\de}(x, a)}{A_{\de}(x, a)}
\eeq
defines an element of the ring $\wh\La_n$. Furthermore,
setting $\wh s_{\la}(x\vt a)=0$ if the length of $\la$
exceeds the number of the $x$ variables, we obtain that
the evaluation of the element $\wh s_{\la}(x\vt a)\in\wh\La_n$ at
$x_n=0$ yields the corresponding element of $\wh\La_{n-1}$
associated with $\la$.
Thus, the sequence $\wh s_{\la}(x\vt a)\in\wh\La_n$
for $n=0,1,\dots$
defines an element $\wh s_{\la}(x\vt a)$ of $\wh\La(x\vt a)$
which we call the {\it dual Schur function\/}.
The lowest degree component of $\wh s_{\la}(x\vt a)$
in $x$ coincides with the Schur function $s_{\la}(x)$.
Moreover, if $a$ is specialized to
the sequence of zeros, then
$\wh s_{\la}(x\vt a)$ specializes to $s_{\la}(x)$.

Now we prove an analogue of the Cauchy identity
involving the double and dual Schur functions.
Consider one more set of variables $y=(y_1,y_2,\dots)$.

\bth\label{thm:cch}
The following identity holds
\beql{cchid}
\prod_{i,\ts j\geqslant 1}\frac{1-a_i\ts y_j}{1-x_i\ts y_j}
=\sum_{\la\in\Pc} s_{\la}(x\vt a)\ts \wh s_{\la}(y\vt a).
\eeq
\eth

\bpf
We use a modification
of the argument applied in \cite[Section~I.4]{m:sfh}
for the proof of the classical Cauchy identity
(see formula (4.3) there).
As we pointed out in Section~\ref{subsec:dpr},
it will be sufficient to prove the identity in the case
of finite sets of variables $x=(x_1,\dots,x_n)$ and
$y=(y_1,\dots,y_n)$. We have
\beql{adexy}
A_{\de}(x\vt a)\ts A_{\de}(y, a)
\sum_{\la\in\Pc} s_{\la}(x\vt a)\ts \wh s_{\la}(y\vt a)
=\sum_{\ga} A_{\ga}(x\vt a)\ts A_{\ga}(y, a),
\eeq
summed over $n$-tuples $\ga=(\ga_1,\dots,\ga_n)$
with $\ga_1>\dots>\ga_n\geqslant 0$. Since
\ben
A_{\ga}(y, a)=\sum_{\si\in\Sym_n}\sgn\tss\si\ts
\prod_{i=1}^n (y_i,a)^{\ga_{\si(i)}}
\tss(1-a^{}_{n-\ga_{\si(i)}-1}y_i)\dots(1-a^{}_{1-\ga_{\si(i)}}y_i)
\een
and $A_{\ga}(x\vt a)$ is skew-symmetric under
permutations of the components of $\ga$, we can write
\eqref{adexy} in the form
\beql{abexa}
\sum_{\be} A_{\be}(x\vt a)
\prod_{i=1}^n (y_i,a)^{\be_i}
\tss(1-a^{}_{n-\be_i-1}y_i)\dots(1-a^{}_{1-\be_i}y_i),
\eeq
summed over $n$-tuples $\be=(\be_1,\dots,\be_n)$
on nonnegative integers. Due to the Jacobi--Trudi
formula \eqref{jt}, we have
\ben
A_{\be}(x\vt a)=A_{\de}(x\vt a)\sum_{\si\in\Sym_n}\sgn\tss\si\cdot
h_{\be_{\si(1)}-n+1}(x\vt a)\dots h_{\be_{\si(n)}}(x\vt \tau^{n-1}a).
\een
Hence, \eqref{abexa} becomes
\begin{multline}
\label{axy}
A_{\de}(x\vt a)
\sum_{\al}h_{\al_1}(x\vt a)\dots h_{\al_n}(x\vt \tau^{n-1}a)\\
{}\times\sum_{\si\in\Sym_n}\sgn\tss\si\cdot
\prod_{i=1}^n (y_{\si(i)},a)^{\al_i+n-i}
\tss(1-a^{}_{i-\al_i-1}y_{\si(i)})
\dots(1-a^{}_{i-\al_i-n+1}y_{\si(i)}),
\end{multline}
summed over $n$-tuples $\al=(\al_1,\dots,\al_n)$
on nonnegative integers. However, using \eqref{genh},
for each $i=1,\dots,n$ we obtain
\ben
\bal
\sum_{k=0}^{\infty}& h_k(x\vt \tau^{i-1}a)\ts (z,a)^{k+n-i}
\tss(1-a^{}_{i-\al_i-1}z)
\dots(1-a^{}_{i-\al_i-n+1}z)\\
{}&=z^{n-i}\ts (1-a_1\tss z)\dots(1-a_{i-1}\tss z)
\sum_{k=0}^{\infty} h_k(x\vt \tau^{i-1}a)\ts (z,\tau^{i-1}a)^k\\
{}&=z^{n-i}\ts (1-a_1\tss z)\dots(1-a_{i-1}\tss z)
\prod_{r=1}^{n}\frac{1-a_{i+r-1}\tss z}{1-x_r\tss z},
\eal
\een
where we put $z=y_{\si(i)}$. Therefore, \eqref{axy} simplifies to
\ben
\bal
A_{\de}(x\vt a)&\ts
\prod_{i,j=1}^{n}\frac{1-a_{i}\tss y_j}{1-x_i\tss y_j}
\sum_{\si\in\Sym_n}\sgn\tss\si\cdot
\prod_{i=1}^n y_{\si(i)}^{n-i}
\ts(1-a^{}_{n+1}y_{\si(i)})
\dots(1-a^{}_{n+i-1}y_{\si(i)})\\
{}&=A_{\de}(x\vt a)\tss A_{\de}(y, a)\ts
\prod_{i,j=1}^{n}\frac{1-a_{i}\tss y_j}{1-x_i\tss y_j},
\eal
\een
thus completing the proof.
\epf

Let $z=(z_1,z_2,\dots)$ be another set of variables.

\bco\label{cor:ssduco}
The following identity holds
\ben
\prod_{i,\ts j\geqslant 1}\frac{1+y_i\ts z_j}{1-x_i\ts z_j}
=\sum_{\la\in\Pc} s_{\la}(x/y\vt a)\ts \wh s_{\la}(z\vt a),
\een
\eco

\bpf
Observe that
the elements $\wh s_{\la}(z\vt a)\in \wh\La(z\vt a)$ are uniquely
determined by this relation. Hence, the claim follows
by the application of Proposition~\ref{prop:imisom}
and Theorem~\ref{thm:cch}.
\epf

Some other identities of this kind are immediate from
the symmetry property \eqref{symprsu}
and Corollary~\ref{cor:ssduco}.

\bco\label{cor:duco}
We have the identities
\ben
\prod_{i,\ts j\geqslant 1}\frac{1+x_i\ts z_j}{1-y_i\ts z_j}
=\sum_{\la\in\Pc} s_{\la}(x/y\vt a)\ts \wh s_{\la'}(z\vt a')
\een
and
\ben
\prod_{i,\ts j\geqslant 1}\frac{1+x_i\ts y_j}{1+a_i\ts y_j}
=\sum_{\la\in\Pc} s_{\la}(x\vt a)\ts \wh s_{\la'}(y\vt a').
\een
\vskip-1.2\baselineskip
\qed
\eco

\subsection{Combinatorial presentation}

Given a skew diagram $\theta$, introduce the
corresponding {\it skew dual Schur function\/}
$\wh s_{\theta}(x\vt a)$ by the formula
\beql{skewd}
\wh s_{\theta}(x\vt a)=\sum_{T}
\prod_{\al\in\theta}X_{T(\al)}(a^{}_{-c(\al)+1},a^{}_{-c(\al)}),
\eeq
summed over the reverse $\theta$-tableaux $T$, where
\ben
X_i(g,h)=
\frac{x_i\ts(1-g\ts x_{i-1})\dots (1-g\ts x_1)}
{(1-h\ts x_{i})\dots (1-h\ts x_1)}.
\een

\bth\label{thm:tab}
For any partition $\mu$ the following identity holds
\beql{skcid}
\prod_{i,\ts j\geqslant 1}\frac{1-a_{i-\mu_i}\ts y_j}{1-x_i\ts y_j}
\ts s_{\mu}(x\vt a)
=\sum_{\nu} s_{\nu}(x\vt a)\ts \wh s_{\nu/\mu}(y\vt a),
\eeq
summed over partitions $\nu$ containing $\mu$.
In particular, if
$\theta=\la$ is a normal {\rm(}nonskew{\rm)} diagram, then
the dual Schur function $\wh s_{\la}(x\vt a)$ admits
the tableau presentation
\eqref{skewd}.
\eth

\bpf
It will be sufficient to consider the case where
the set of variables $y$ is finite,
$y=(y_1,\dots,y_n)$. We will argue by induction on $n$
and suppose that $n\geqslant 1$. By the induction hypothesis,
the identity \eqref{skcid} holds
for the set of variables $y'=(y_2,\dots,y_n)$.
Hence, we need to verify that
\ben
\prod_{i\geqslant 1}\frac{1-a_{i-\mu_i}\ts y_1}{1-x_i\ts y_1}
\ts\sum_{\la} s_{\la}(x\vt a)\ts \wh s_{\la/\mu}(y'\vt a)
=\sum_{\nu} s_{\nu}(x\vt a)\ts \wh s_{\nu/\mu}(y\vt a).
\een
However, due to \eqref{genh},
\ben
\prod_{i\geqslant 1}\frac{1-a_{i-\mu_i}\ts y_1}{1-x_i\ts y_1}
=\sum_{k=0}^{\infty} \frac{h_k(x\vt a^{\mu})\ts y_1^k}
{(1-a_0\tss y_1)\dots(1-a_{-k+1}\tss y_1)},
\een
where $a^{\mu}$ denotes the sequence of parameters such that
$(a^{\mu})_i=a_{i-\mu_i}$ for $i\geqslant 1$ and
$(a^{\mu})_i=a_i$ for $i\leqslant 0$. Now define
polynomials $c_{\la,(k)}^{\nu}(a,a^{\mu})\in\QQ[a]$ by
the expansion
\ben
s_{\la}(x\vt a)\ts h_k(x\vt a^{\mu})=\sum_{\nu}
c_{\la,(k)}^{\nu}(a,a^{\mu})\ts s_{\nu}(x\vt a).
\een
Hence, the claim will follow if we show that
\beql{snumu}
\wh s_{\nu/\mu}(y\vt a)=\sum_{\la,\ts k}c_{\la,(k)}^{\nu}(a,a^{\mu})
\ts \wh s_{\la/\mu}(y'\vt a)\ts \frac{y_1^k}
{(1-a_0\tss y_1)\dots(1-a_{-k+1}\tss y_1)}.
\eeq
The definition \eqref{skewd} of the skew dual Schur
functions implies that
\ben
\wh s_{\nu/\mu}(y\vt a)=\sum_{\la}
\wh s_{\la/\mu}(y'\vt a)\ts\prod_{\al\in\la/\mu}
\frac{1-a_{-c(\al)+1}\tss y_1}{1-a_{-c(\al)}\tss y_1}
\ts\prod_{\be\in\nu/\la}
\frac{y_1}{1-a_{-c(\al)}\tss y_1},
\een
summed over diagrams $\la$ such that
$\mu\subseteq\la\subseteq\nu$ and $\nu/\la$ is a horizontal
strip (i.e., every column of this diagram contains
at most one box). Therefore, \eqref{snumu}
will follow from the relation
\ben
\bal
\sum_{k}&c_{\la,(k)}^{\nu}(a,a^{\mu})
\ts\frac{y_1^k}
{(1-a_0\tss y_1)\dots(1-a_{-k+1}\tss y_1)}\\
{}&=\prod_{\al\in\la/\mu}
\frac{1-a_{-c(\al)+1}\tss y_1}{1-a_{-c(\al)}\tss y_1}
\ts\prod_{\be\in\nu/\la}
\frac{y_1}{1-a_{-c(\al)}\tss y_1}
\eal
\een
which takes more convenient form after the
substitution $t=y_1^{-1}$:
\beql{sumts}
\sum_{k}
\ts\frac{c_{\la,(k)}^{\nu}(a,a^{\mu})}
{(t-a_0)\dots(t-a_{-k+1})}=\prod_{\al\in\la/\mu}
\big(t-a_{-c(\al)+1}\big)
\ts\prod_{\be\in\nu/\mu}
\big(t-a_{-c(\al)}\big)^{-1}.
\eeq
We will verify the latter by induction on $|\nu|-|\la|$.
Suppose first that $\nu=\la$. Then
$c_{\la,(k)}^{\la}(a,a^{\mu})=h_k(a_{\la}\vt a^{\mu})$
by Proposition~\ref{prop:interp}, and relation \eqref{genh}
implies that
\ben
\sum_{k}
\ts\frac{h_k(a_{\la}\vt a^{\mu})}
{(t-a_0)\dots(t-a_{-k+1})}=\prod_{i\geqslant 1}
\frac{t-a_{i-\mu_i}}{t-a_{i-\la_i}}.
\een
This expression coincides with
\ben
\prod_{\al\in\la/\mu}
\frac{t-a_{-c(\al)+1}}{t-a_{-c(\al)}},
\een
thus verifying \eqref{sumts} in the case under
consideration. Suppose now that $|\nu|-|\la|\geqslant 1$.
By Proposition~\ref{prop:interp},
we have
\ben
c_{\la,(k)}^{\nu}(a,a^{\mu})=\frac{1}{|a_{\nu}|-|a_{\la}|}\ts
\Bigg(\sum_{\la^+,\ts\la\to\la^+}c_{\la^+,(k)}^{\nu}(a,a^{\mu})
-\sum_{\nu^-,\ts\nu^-\to\nu}c_{\la,(k)}^{\nu^-}(a,a^{\mu})\Bigg).
\een
Hence, applying the induction hypothesis, we can
write the left hand side of \eqref{sumts} in the form
\ben
\bal
\frac{1}{|a_{\nu}|-|a_{\la}|}\ts
\Big(&\sum_{\la^+}\prod_{\al\in\la^+/\mu}
\big(t-a_{-c(\al)+1}\big)
\ts\prod_{\be\in\nu/\mu}
\big(t-a_{-c(\al)}\big)^{-1}\\
{}-{}&\sum_{\nu^-}\prod_{\al\in\la/\mu}
\big(t-a_{-c(\al)+1}\big)
\ts\prod_{\be\in\nu^-/\mu}
\big(t-a_{-c(\al)}\big)^{-1}\Big).
\eal
\een
Since $\nu/\la$ is a horizontal strip, we have
\ben
\sum_{\al=\la^+/\la}\big(t-a_{-c(\al)+1}\big)-
\sum_{\al=\nu/\nu^-}\big(t-a_{-c(\al)}\big)=|a_{\nu}|-|a_{\la}|,
\een
so that the previous expression simplifies to
\ben
\prod_{\al\in\la/\mu}
\big(t-a_{-c(\al)+1}\big)
\ts\prod_{\be\in\nu/\mu}
\big(t-a_{-c(\al)}\big)^{-1}
\een
completing the proof of \eqref{sumts}.

The second part of the proposition follows from
Theorem~\ref{thm:cch} and the fact that the
elements $\wh s_{\la}(y\vt a)\in\wh\La(y\vt a)$ are uniquely
determined by the relation \eqref{cchid}.
\epf

\bre\label{rem:aze}
Under the specialization $a_i=0$ the identity
of Theorem~\ref{thm:tab} turns into a particular case
of the identity in \cite[Example~I.5.26]{m:sfh}.
\qed
\ere

Since the skew dual Schur functions are uniquely determined
by the expansion \eqref{skcid}, the following corollary
is immediate from Theorem~\ref{thm:tab}.

\bco\label{cor:symme}
The skew dual Schur functions defined in \eqref{skewd}
belong to the ring $\wh\La(x\vt a)$. In particular,
they are symmetric in the variables $x$.
\qed
\eco

Recall the Littlewood--Richardson polynomials
defined by \eqref{lrpoldef}.

\bpr\label{prop:skelrp}
For any skew diagram $\nu/\mu$ we have the expansion
\ben
\wh s_{\nu/\mu}(y\vt a)=\sum_{\la} c_{\la\mu}^{\tss\nu}(a)\ts
\wh s_{\la}(y\vt a).
\een
\epr

\bpf
We use an argument similar to the one used
in \cite[Section~I.5]{m:sfh}.
Consider the set of variables $(y,y^{\tss\prime})$,
where $y=(y_1,y_2,\dots)$
and $y^{\tss\prime}=(y^{\tss\prime}_1,y^{\tss\prime}_2,\dots)$
and assume they are ordered in the way
that each $y_i$ precedes each $y^{\tss\prime}_j$. By the tableau
presentation \eqref{skewd} of the dual Schur functions,
we get
\beql{snuyza}
\wh s_{\nu}(y,y^{\tss\prime}\vt a)=\sum_{\mu\subseteq\nu}
\wh s_{\nu/\mu}(y\vt a)\ts \wh s_{\mu}(y^{\tss\prime}\vt a).
\eeq
On the other hand, by Theorem~\ref{thm:cch},
\ben
\bal
{}&\sum_{\nu}s_{\nu}(x\vt a)\ts \wh s_{\nu}(y,y^{\tss\prime}\vt a)=
\prod_{i,j\geqslant 1}\frac{1-a_i\tss y_j}{1-x_i\tss y_j}
\ts\prod_{i,k\geqslant 1}
\frac{1-a_i\tss y^{\tss\prime}_k}{1-x_i\tss y^{\tss\prime}_k}\\
{}&=\sum_{\la,\ts\mu}s_{\la}(x\vt a)\ts \wh s_{\la}(y\vt a)\ts
s_{\mu}(x\vt a)\ts \wh s_{\mu}(y^{\tss\prime}\vt a)
=\sum_{\la,\ts\mu,\ts\nu}c_{\la\mu}^{\tss\nu}(a)\ts
s_{\nu}(x\vt a)\ts \wh s_{\la}(y\vt a)\ts
\wh s_{\mu}(y^{\tss\prime}\vt a)
\eal
\een
which proves that
\beql{snulrp}
\wh s_{\nu}(y,y^{\tss\prime}\vt a)=
\sum_{\la,\ts\mu}c_{\la\mu}^{\tss\nu}(a)\ts
\wh s_{\la}(y\vt a)\ts \wh s_{\mu}(y^{\tss\prime}\vt a).
\eeq
The desired relation now follows by comparing \eqref{snuyza}
and \eqref{snulrp}.
\epf

\subsection{Jacobi--Trudi-type formulas}

Introduce the {\it dual elementary\/} and {\it complete
symmetric functions\/} by
\ben
\wh e_k(x\vt a)=\wh s_{(1^k)}(x\vt a),\qquad
\wh h_k(x\vt a)=\wh s_{(k)}(x\vt a).
\een
By Theorem~\ref{thm:tab},
\ben
\bal
\wh e_k(x\vt a)&=\sum_{i_1>\dots>i_k}
X_{i_1}(a_1,a_0)X_{i_2}(a_2,a_1)\dots
X_{i_k}(a_k,a_{k-1}),\\
\wh h_k(x\vt a)&=\sum_{i_1\geqslant\dots\geqslant i_k}
X_{i_1}(a_1,a_0)X_{i_2}(a_0,a_{-1})\dots
X_{i_k}(a_{-k+2},a_{-k+1}).
\eal
\een

\bpr\label{prop:gendu}
We have the following generating series formulas
\ben
\bal
1+\sum_{k=1}^{\infty}\wh e_k(x\vt a)\ts
(t+a_0)(t+a_1)\dots (t+a_{k-1})=
\prod_{i=1}^{\infty}\frac{1+t\tss x_i}{1-a_0\tss x_i},\\
1+\sum_{k=1}^{\infty}\wh h_k(x\vt a)\ts
(t-a_1)(t-a_0)\dots (t-a_{-k+2})=
\prod_{i=1}^{\infty}\frac{1-a_1\tss x_i}{1-t\tss x_i}.
\eal
\een
\epr

\bpf
The first relation follows from
the second identity in Corollary~\ref{cor:duco}
by taking $x=(t)$ and then replacing $a$ by $a'$ and
$y_i$ by $x_i$ for all $i$. Similarly,
the second relation follows from Theorem~\ref{thm:cch}
by taking $x=(t)$ and replacing
$y_i$ by $x_i$.
\epf

We can now prove an analogue of the Jacobi--Trudi formula
for the dual Schur functions.

\bpr\label{prop:jtdu}
If $\la$ and $\mu$ are
partitions of length not exceeding $n$, then
\beql{jtdu}
\wh s_{\la/\mu}(x\vt a)
=\det\big[\tss\wh h_{\la_i-\mu_j-i+j}
(x\vt \tau^{\tss -\mu_j+j-1}\tss a)\big]_{i,j=1}^n.
\eeq
\epr

\bpf
Apply Theorem~\ref{thm:tab} for the finite set
of variables $x=(x_1,\dots,x_n)$ and multiply both
sides of \eqref{skcid} by $A_{\de}(x\vt a)$.
This gives
\beql{pramu}
\prod_{j\geqslant 1}\prod_{i=1}^n
\frac{1-a_{i-\mu_i}\ts y_j}{1-x_i\ts y_j}
\ts A_{\mu+\de}(x\vt a)
=\sum_{\la} A_{\la+\de}(x\vt a)\ts \wh s_{\la/\mu}(y\vt a).
\eeq
For any $\si\in\Sym_n$ we have
\ben
\prod_{j\geqslant 1}\prod_{i=1}^n
\frac{1-a_{i-\mu_i}\ts y_j}{1-x_i\ts y_j}
=\prod_{j\geqslant 1}\prod_{i=1}^n
\frac{1-a_{i-\mu_i}\ts y_j}{1-x_{\si(i)}\ts y_j}.
\een
By the second formula of Proposition~\ref{prop:gendu},
\ben
\prod_{j\geqslant 1}
\frac{1-a_{i-\mu_i}\ts y_j}{1-x_{\si(i)}\ts y_j}
=\sum_{k=0}^{\infty}
\wh h_k(y\vt\tau^{-\mu_i+i-1}a)\ts
(x_{\si(i)}-a_{i-\mu_i})\dots(x_{\si(i)}-a_{i-\mu_i-k+1}).
\een
Since
\ben
A_{\mu+\de}(x\vt a)=\sum_{\si\in\Sym_n}\sgn\si\cdot
(x_{\si(1)}\vt a)^{\mu_1+n-1}\dots (x_{\si(n)}\vt a)^{\mu_n} ,
\een
the left hand side of \eqref{pramu}
can be written in the form
\ben
\sum_{\si\in\Sym_n}\sgn\si\ts
\prod_{i=1}^n\sum_{k_i=0}^{\infty}
(x_{\si(i)}-a_n)\dots(x_{\si(i)}-a_{i-\mu_i-k_i+1})
\ts h_{k_i}(y,\tau^{-\mu_i+i-1}a).
\een
Hence, comparing the coefficients
of $(x_{1}\vt a)^{\la_1+n-1}\dots (x_n\vt a)^{\la_n}$
on both sides of \eqref{pramu}, we get
\ben
\wh s_{\la/\mu}(y\vt a)=\sum_{\rho\in\Sym_n}
\sgn\rho\ts\prod_{i=1}^n \wh h_{\la_i-\mu_{\rho(i)}-i+\rho(i)}
(y\vt\tau^{-\mu_{\rho(i)}+\rho(i)-1}a),
\een
as required.
\epf

Proposition~\ref{prop:jtdu} implies that the
dual Schur functions may be regarded as a specialization
of the generalized Schur functions described in
\cite[9th~Variation]{m:sf}. Namely,
in the notation of that paper, specialize the variables
$h_{rs}$ by
\beql{hrsspe}
h_{rs}=\wh h_r(x\vt\tau^{-s}a),\qquad r\geqslant 1,\quad s\in\ZZ.
\eeq
Then the Schur functions $s_{\la/\mu}$ of \cite{m:sf}
become $\wh s_{\la/\mu}(x\vt a)$.
Hence the following corollaries are immediate
from $(9.6')$ and $(9.7)$ in \cite{m:sf}
and Proposition~\ref{prop:jtdu}.
The first of them is an analogue of the
N\"{a}gelsbach--Kostka formula.

\bco\label{cor:nkdu}
If $\la$ and $\mu$ are
partitions such that the lengths
of $\la'$ and $\mu'$ do not exceed $m$, then
\beql{nkdu}
\wh s_{\la/\mu}(x\vt a)
=\det\big[\wh e_{\la'_i-\mu'_j-i+j}
(x\vt\tau^{\tss \mu'_j-j+1}\tss a)\big]_{i,j=1}^m.
\eeq
\eco

Suppose that $\la$ is a diagram with $d$ boxes
on the main diagonal. Write $\la$ in the
Frobenius notation
\ben
\la=(\al_1,\dots,\al_d\tss|\tss\be_1,\dots,\be_d)=(\al\tss|\tss\be),
\een
where $\al_i=\la_i-i$ and $\be_i=\la'_i-i$.
The following is an analogue of the Giambelli formula.

\bco\label{cor:gidu}
We have the identity
\beql{gidu}
\wh s_{(\al\tss|\tss\be)}(x\vt a)
=\det\big[\wh s_{(\al_i\tss|\tss\be_j)}
(x\vt a)\big]_{i,j=1}^d.
\eeq
\eco

\subsection{Expansions in terms of Schur functions}

We will now deduce
expansions of the dual Schur functions in terms
of the Schur functions $s_{\la}(x)$
whose coefficients are elements of $\QQ[a]$
written explicitly as certain determinants.
In Theorem~\ref{thm:tabre} below we will give
alternative tableau presentations for these coefficients.
Suppose that $\mu$ is a diagram containing $d$ boxes
on the main diagonal.

\bpr\label{prop:dsfexp}
The dual Schur function $\wh s_{\mu}(x\vt a)$ can be written
as the series
\ben
\bal
\wh s_{\mu}(x\vt a)=\sum_{\la}(-1)^{n(\la/\mu)}\ts
&\det\big[h_{\la_i-\mu_j-i+j}
(a^{}_0,a^{}_{-1},\dots,a^{}_{j-\mu_j})\big]_{i,j=1}^d\\
{}\times{}&
\det\big[e_{\la_i-\mu_j-i+j}
(a^{}_1,a^{}_2,\dots,a^{}_{j-\mu_j-1})\big]_{i,j=d+1}^n
\ts s_{\la}(x),
\eal
\een
summed over diagrams $\la$ which
contain $\mu$ and such that $\la$ has
$d$ boxes on the main diagonal, where
$n(\la/\mu)$ denotes the total number of boxes in
the diagram $\la/\mu$ in rows $d+1,d+2,\dots,n=\ell(\la)$.
\epr

\bpf
It will be sufficient to prove the formula for the case
of finite set of variables $x=(x_1,\dots,x_n)$.
We use the definition \eqref{dusf} of the dual Schur
functions. The entries $A_{ij}$ of
the determinant $A_{\mu+\de}(x, a)$
can be written as
\ben
A_{ij}=\begin{cases}{\displaystyle\sum_{p_j\geqslant 0}
h_{p_j}(a^{}_0,a^{}_{-1},\dots,a^{}_{j-\mu_j})\ts
x_i^{\mu_j+p_j+n-j}}
\quad&\text{for}\quad j=1,\dots,d,\\[2em]
{\displaystyle\sum_{p_j\geqslant 0}
(-1)^{p_j}\ts
e_{p_j}(a^{}_1,a^{}_2,\dots,a^{}_{j-\mu_j-1})\ts
x_i^{\mu_j+p_j+n-j}}
\quad&\text{for}\quad j=d+1,\dots,n.
\end{cases}
\een
Hence, \eqref{dusf} gives
\ben
\bal
\wh s_{\mu}(x\vt a)={}&\sum_{p_1,\dots,\ts p_n}
\prod_{j=1}^d
h_{p_j}(a^{}_0,a^{}_{-1},\dots,a^{}_{j-\mu_j})
\prod_{j=d+1}^n
(-1)^{p_j}\ts e_{p_j}(a^{}_1,a^{}_2,\dots,a^{}_{j-\mu_j-1})\\[1em]
{}&{}\times\det[x_i^{\mu_j+p_j+n-j}]/\det[x_i^{n-j}].
\eal
\een
The ratio of the determinants in this formula
is nonzero only if
\ben
\mu_{\si(j)}+p_{\si(j)}+n-{\si(j)}=\la_j+n-j,\qquad j=1,\dots,n,
\een
for some diagram $\la$ containing $\mu$ and some
permutation $\si$ of the set $\{1,\dots,n\}$.
Moreover, since $e_{p_j}(a^{}_1,a^{}_2,\dots,a^{}_{j-\mu_j-1})=0$
for $p_j>j-\mu_j-1$, the number of diagonal boxes
in $\la$ equals $d$.
The ratio can then be written as
\ben
\det[x_i^{\mu_j+p_j+n-j}]/\det[x_i^{n-j}]=\sgn\si\cdot
s_{\la}(x),
\een
which gives the desired formula for the coefficients.
\epf

\bco\label{cor:sone}
Using the Frobenius notation $(\al\tss|\tss\be)$
for the hook diagram $(\al+1,1^{\be})$, we
have
\ben
\wh s_{(\al\tss|\tss\be)}(x\vt a)=
\sum_{p,\ts q\geqslant 0}(-1)^q\ts
h_p(a^{}_0,a^{}_{-1},\dots,a^{}_{-\al})\ts
h_q(a^{}_1,a^{}_2,\dots,a^{}_{\be+1})\ts
s^{}_{(\al+p\tss|\tss\be+q)}(x).
\een
\eco

\bpf
By Proposition~\ref{prop:dsfexp}, the coefficient of
$s^{}_{(\al+p\tss|\tss\be+q)}(x)$ in the expansion
of the dual Schur function $\wh s_{(\al\tss|\tss\be)}(x\vt a)$ equals
\beql{hdete}
(-1)^q\ts
h_p(a^{}_0,a^{}_{-1},\dots,a^{}_{-\al})\ts
\det\big[e_{j-i+1}
(a^{}_1,a^{}_2,\dots,a^{}_{\be+j})\big]_{i,j=1}^q.
\eeq
Using the relations for the elementary symmetric
polynomials
\ben
e_k(a^{}_1,a^{}_2,\dots,a^{}_{\be+j})
=e_k(a^{}_1,a^{}_2,\dots,a^{}_{\be+j-1})+
e_{k-1}(a^{}_1,a^{}_2,\dots,a^{}_{\be+j-1})\ts a^{}_{\be+j},
\een
it is not difficult to bring the determinant
which occurs in \eqref{hdete} to the form
\beql{deteq}
\det\big[e_{j-i+1}
(a^{}_1,a^{}_2,\dots,a^{}_{\be+j})\big]_{i,j=1}^q
=\det\big[e_{j-i+1}
(a^{}_1,a^{}_2,\dots,a^{}_{\be+1})\big]_{i,j=1}^q.
\eeq
Indeed, denote by $C_1,\dots,C_q$ the columns of the
$q\times q$ matrix which occurs
on the left hand side. Now replace
$C_j$ by $C_j-a^{}_{\be+j}\ts C_{j-1}$
consequently for $j=q,q-1,\dots,2$. These operations leave
the determinant of the matrix unchanged, while for $j\geqslant 2$ the
$(i,j)$ entry of the new matrix equals
$e_{j-i+1}
(a^{}_1,a^{}_2,\dots,a^{}_{\be+j-1})$. Applying similar
column operations to the new matrix and using obvious induction
we will bring its determinant
to the form which occurs on the right hand side of \eqref{deteq}.
However, this determinant
coincides with $h_q(a^{}_1,a^{}_2,\dots,a^{}_{\be+1})$
due to the N\"{a}gelsbach--Kostka formula (i.e., \eqref{nkdu}
with the zero sequence $a$; that is, $a_i=0$ for all $i\in\ZZ$).
\epf

\bex\label{ex:sone}
The dual Schur function corresponding
to the single box diagram is given by
\ben
\wh s_{(1)}(x\vt a)=\sum_{p,\ts q\geqslant 0}
(-1)^q\ts a_0^p\ts a_1^q\ts
s^{}_{(p\tss|\tss q)}(x).
\een
\vskip-1.2\baselineskip
\qed
\eex

Recall that the involution
$\om:\La\to\La$ on the ring of symmetric functions
in $x$ takes $s_{\la}(x)$ to $s_{\la'}(x)$; see
\cite[Section~I.2]{m:sfh} or Section~\ref{sec:def} above.
Let us extend $\om$ to the $\QQ[a]$-linear
involution
\beql{whome}
\wh\om:\wh\La(x\vt a)\to \wh\La(x\vt a),\qquad
\sum_{\la\in\Pc} c_{\la}(a)\ts s_{\la}(x)\mapsto
\sum_{\la\in\Pc} c_{\la}(a)\ts s_{\la'}(x),
\eeq
where $c_{\la}(a)\in\QQ[a]$.
We will find the images
of the dual Schur functions under $\wh\om$.
As before, by $a'$ we denote the sequence of variables
such that $(a')_i=-a_{-i+1}$ for all $i\in\ZZ$.

\bco\label{cor:omhat}
For any skew diagram $\la/\mu$ we have
\beql{imomdu}
\wh\om:\wh s_{\la/\mu}(x\vt a)\mapsto  \wh s_{\la'/\mu'}(x\vt a').
\eeq
\eco

\bpf
By Corollary~\ref{cor:sone}, for any $m\in\ZZ$
\ben
\wh\om:
\wh s_{(\al\tss|\tss\be)}(x\vt\tau^ma)\mapsto
\wh s_{(\be\tss|\tss\al)}(x\vt\tau^{-m}a').
\een
In particular,
\ben
\wh\om:
\wh h_k(x\vt\tau^m a)\mapsto \wh e_k(x\vt\tau^{-m}a'),
\qquad k\geqslant 0.
\een
The statement now follows from \eqref{jtdu} and \eqref{nkdu}.
\epf

Note that \eqref{imomdu} with $\mu=\varnothing$ also
follows from Corollary~\ref{cor:sone}
and the Giambelli formula \eqref{gidu}.

\bre\label{rem:athom}
The involution $\wh\om$ does not coincide
with the involution introduced in \cite[(9.6)]{m:sfh}.
The latter is defined on the ring generated by the
elements $h_{rs}$ and
takes the generalized Schur function
$s_{\la/\mu}$ to $s_{\la'/\mu'}$. Therefore, under the
specialization \eqref{hrsspe}, the image of
$\wh s_{\la/\mu}(x\vt a)$ would be $\wh s_{\la'/\mu'}(x\vt a)$
which is different from \eqref{imomdu}.
\qed
\ere

We can now derive an alternative expansion
of the dual Schur functions in terms of the Schur functions
$s_{\la}(x)$; cf. Proposition~\ref{prop:dsfexp}.
Suppose that $\la$
is a diagram
which contains $\mu$ and such that $\mu$ and $\la$ have the
same number of boxes $d$ on the diagonal.
By a {\it hook $\la/\mu$-tableau\/}
$T$ we will mean a tableau obtained by filling in the boxes
of $\la/\mu$ with integers in the following way.
The entries in the first $d$ rows
weakly increase along the rows and strictly increase
down the columns, and all entries in row $i$
belong to the set $\{i-\mu_i,\dots,-1,0\}$ for
$i=1,\dots,d$;
the entries in the first $d$ columns
weakly decrease down the columns and strictly decrease
along the rows, and all entries in column $j$
belong to the set $\{1,2,\dots,\mu'_j-j+1\}$ for
$j=1,\dots,d$. Then we define the corresponding flagged
Schur function $\vp_{\la/\mu}(a)$ by the formula
\ben
\vp_{\la/\mu}(a)=\sum_T \prod_{\al\in\la/\mu} a^{}_{T(\al)},
\een
summed over the hook $\la/\mu$-tableaux $T$.

\bth\label{thm:tabre}
Let $\mu$ be a diagram and let $d$ be the number
of boxes on the main diagonal of $\mu$.
We have the expansion of the dual Schur function
$\wh s_{\mu}(x\vt a)$
\ben
\wh s_{\mu}(x\vt a)=\sum_{\la}(-1)^{n(\la/\mu)}\ts
\vp_{\la/\mu}(a)\ts s_{\la}(x),
\een
summed over diagrams $\la$ which
contain $\mu$ and such that $\la$ has
$d$ boxes on the main diagonal, where
$n(\la/\mu)$ denotes the total number of boxes in
the diagram $\la/\mu$ in rows $d+1,d+2,\dots$.
\eth

\bpf
Consider the expansions of $\wh s_{\mu}(x\vt a)$ and
$\wh s_{\mu'}(x\vt a')$
provided by Proposition~\ref{prop:dsfexp}.
By Corollary~\ref{cor:omhat}, $\wh s_{\mu}(x\vt a)$ is the image
of $\wh s_{\mu'}(x\vt a')$ under the involution $\wh\om$.
Since $\wh\om:s_{\la}(x)\mapsto s_{\la'}(x)$,
taking $\la_i=\mu_i$ for $i=1,\dots,d$ and comparing
the coefficients of $s_{\la}(x)$ in
the expansions of $\wh s_{\mu}(x\vt a)$ and
$\wh\om\big(\wh s_{\mu'}(x\vt a')\big)$,
we can conclude that
\begin{multline}
(-1)^{n(\la/\mu)}
\det\big[e_{\la_i-\mu_j-i+j}
(a^{}_1,a^{}_2,\dots,a^{}_{j-\mu_j-1})\big]_{i,j\geqslant d+1}\\
=\det\big[h_{\la'_i-\mu'_j-i+j}
(a'_0,a'_{-1},\dots,a'_{j-\mu'_j})\big]_{i,j=1}^d
\non
\end{multline}
so that
\ben
\det\big[e_{\la_i-\mu_j-i+j}
(a^{}_1,a^{}_2,\dots,a^{}_{j-\mu_j-1})\big]_{i,j\geqslant d+1}
=\det\big[h_{\la'_i-\mu'_j-i+j}
(a^{}_1,a^{}_2,\dots,a^{}_{\mu'_j-j+1})\big]_{i,j=1}^d.
\een
On the other hand, if $\la$ is a diagram containing
$\mu$ and such that $\la$ has
$d$ boxes on the main diagonal,
both determinants
\ben
\det\big[h_{\la_i-\mu_j-i+j}
(a^{}_0,a^{}_{-1},\dots,a^{}_{j-\mu_j})\big]_{i,j=1}^d,\quad
\det\big[h_{\la'_i-\mu'_j-i+j}
(a^{}_1,a^{}_2,\dots,a^{}_{\mu'_j-j+1})\big]_{i,j=1}^d
\een
coincide with
the respective `row-flagged Schur functions'
of \cite[(8.2)]{m:sf},
\cite{w:fs}, and they admit
the required tableau presentations.
\epf

It is clear from the definition of the
flagged Schur function $\vp_{\la/\mu}(a)$ that it
can be written as the product of two polynomials.
More precisely, suppose that
the diagram $\la$ contains
$\mu$ and both $\la$ and $\mu$ have
$d$ boxes on their main diagonals.
Let $(\la/\mu)_+$ denote the part of the skew diagram
$\la/\mu$ contained in the top $d$ rows.
With this notation,
the hook flagged Schur function $\vp_{\la/\mu}(a)$
can be written as
\beql{whsufa}
\vp_{\la/\mu}(a)=(-1)^{n(\la/\mu)}\ts
\vp_{(\la/\mu)_+}(a)\ts\vp_{(\la'/\mu')_+}(a').
\eeq
In addition to the tableau presentation of
the polynomial $\vp_{(\la/\mu)_+}(a)$ given above,
we can get an alternative presentation based on the
column-flagged Schur functions; see \cite[$(8.2')$]{m:sf},
\cite{w:fs}. Due to \eqref{whsufa}, this also gives
alternative formulas for the
coefficients in the expansion of $\wh s_{\mu}(x\vt a)$.

\bco\label{cor:anexpa}
We have the tableau presentation
\ben
\vp_{(\la/\mu)_+}(a)=
\sum_T \prod_{\al\in(\la/\mu)_+} a^{}_{T(\al)},
\een
summed over the $(\la/\mu)_+$-tableaux $T$
whose entries in column $j$ belong to the set
$\{0,-1,\dots,-j+\mu'_j+2\}$ for $j\geqslant d+1$,
and the entries weakly increase along
the rows and strictly increase down the columns.
\qed
\eco

Our next goal is to derive the inverse formulas
expressing the Schur functions $s_{\mu}(x)$
as a series of the dual Schur functions $\wh s_{\la}(x\vt a)$.

\bpr\label{prop:inve}
We have the expansion
\ben
\bal
s_{\mu}(x)=\sum_{\la}(-1)^{m(\la/\mu)}\ts
&\det\big[e_{\la_i-\mu_j-i+j}
(a^{}_0,a^{}_{-1},\dots,a^{}_{i-\la_i+1})\big]_{i,j=1}^d\\
{}\times{}&
\det\big[h_{\la_i-\mu_j-i+j}
(a^{}_1,a^{}_2,\dots,a^{}_{i-\la_i})\big]_{i,j\geqslant d+1}
\ts \wh s_{\la}(x\vt a),
\eal
\een
summed over diagrams $\la$ which
contain $\mu$ and such that $\la$ has
$d$ boxes on the main diagonal,
where $m(\la/\mu)$ denotes the total number of boxes in
the diagram $\la/\mu$ in rows $1,\dots,d$.
\epr

\bpf
We will work with a finite set of
variables $x=(x_1,\dots,x_n)$.
The one variable specialization of \eqref{genh} gives
\ben
1+\sum_{k=1}^{\infty} \frac{(x-a_1)(x-a_0)\dots (x-a_{-k+2})\ts t^k}
{(1-a_0\tss t)\dots(1-a_{-k+1}\tss t)}=
\frac{1-a_1\tss t}{1-x\tss t}.
\een
This implies
\beql{decxt}
\sum_{k=1}^{\infty} \frac{(x-a_0)\dots (x-a_{-k+2})\ts t^k}
{(1-a_0\tss t)\dots(1-a_{-k+1}\tss t)}=
\frac{t}{1-x\tss t}.
\eeq
Writing
\ben
(x-a_0)\dots (x-a_{-k+2})=\sum_{i=1}^k(-1)^{k-i}
e_{k-i}(a_0,a_{-1},\dots,a_{-k+2})\ts x^{i-1}
\een
and comparing the coefficients of $x^{r-1}$ on both sides
of \eqref{decxt} we come to the relation
\beql{tpor}
t^{\tss r}=\sum_{k=r}^{\infty}
\frac{(-1)^{k-r}\ts e_{k-r}(a_0,a_{-1},\dots,a_{-k+2})\ts t^k}
{(1-a_0\tss t)\dots(1-a_{-k+1}\tss t)},\qquad r\geqslant 1.
\eeq
Similarly, writing
\ben
\frac{1}
{(1-a_0\tss t)\dots(1-a_{-k+1}\tss t)}
=\sum_{j=0}^{\infty}
h_{j}(a_0,a_{-1},\dots,a_{-k+1})\ts t^{j}
\een
and comparing the coefficients of $t^{\tss r+1}$ on both sides
of \eqref{decxt} we come to
\beql{xpor}
x^r=\sum_{k=0}^r h_{r-k}(a_0,a_{-1},\dots,a_{-k})
\tss(x-a_0)(x-a_{-1})\dots(x-a_{-k+1}),\qquad r\geqslant 0.
\eeq
Assuming that the length of $\mu$ does not exceed $n$,
represent $s_{\mu}(x)$ as the ratio of determinants
\ben
s_{\mu}(x)=\frac{A_{\mu+\de}(x)}{A_{\de}(x)},
\een
where
\ben
A_{\al}(x)=\det\big[x_i^{\al_j}\big]_{i,j=1}^n,
\qquad \al=(\al_1,\dots,\al_n).
\een
By \eqref{tpor},
for any $j=1,\dots,d$ we have
\ben
x_i^{\mu_j-j+1}=\sum_{p=\mu_j-j+1}^{\infty}
\frac{(-1)^{p-\mu_j+j-1}\ts e_{p-\mu_j+j-1}
(a_0,a_{-1},\dots,a_{-p+2})\ts x_i^p}
{(1-a_0\tss x_i)\dots(1-a_{-p+1}\tss x_i)}.
\een
Similarly, for $j=d+1,\dots,n$ we find from \eqref{xpor}
applied for $x=x_i^{-1}$ and $r=j-\mu_j-1$ that
\ben
x_i^{\mu_j-j+1}=\sum_{p=0}^{j-\mu_j-1}
h_{j-\mu_j-p-1}(a_1,a_2,\dots,a_{p+1})\ts x_i^{-p}
\tss(1-a_1 x_i)(1-a_2 x_i)\dots(1-a_p x_i).
\een
Multiplying both sides of these relations by $x_i^{n-1}$
we get the respective expansions of $x_i^{\mu_j+n-j}$
which allow us to write
\ben
\bal
A_{\mu+\de}(x)=\sum_{\be_1,\dots,\ts \be_n}
{}&{}\prod_{j=1}^d
(-1)^{\be_j-\mu_j+j-1}\ts
e_{\be_j-\mu_j+j-1}(a^{}_0,a^{}_{-1},\dots,a^{}_{-\be_j+2})\\
{}\times{}&{}\prod_{j=d+1}^n
\ts h_{j-\mu_j-\be_j-1}(a^{}_1,a^{}_2,\dots,a^{}_{\be_j+1})
\ts A_{\be}(x, a).
\eal
\een
Nonzero summands here correspond to the $n$-tuples $\be$
of the form
\ben
\be_j=\la_{\si(j)}-\si(j)+1,\qquad j=1,\dots,d,
\een
and
\ben
\be_j=-\la_{\tau(j)}+\tau(j)-1,\qquad j=d+1,\dots,n,
\een
where $\si$ is a permutation of $\{1,\dots,d\}$ and $\tau$
is a permutation of
$\{d+1,\dots,n\}$, and
$\la$ is a diagram containing $\mu$
such that $\la$ has $d$ boxes on the main diagonal.
Dividing both sides of the above relation by the Vandermonde
determinant, we get the desired expansion formula.
\epf

Now we obtain a tableau presentation
of the coefficients in the expansion of $s_{\mu}(x)$;
cf. Theorem~\ref{thm:tabre}. We assume, as before, that
$\la$ and $\mu$ have the same number of boxes $d$ on their
main diagonals.
By a {\it dual hook $\la/\mu$-tableau\/}
$T$ we will mean a tableau obtained by filling in the boxes
of $\la/\mu$ with integers in the following way.
The entries in the first $d$ rows
strictly decrease along the rows and weakly decrease
down the columns, and all entries in row $i$
belong to the set $\{0,-1,\dots,i-\la_i+1\}$ for $i=1,\dots,d$;
the entries in the first $d$ columns
strictly increase down the columns and weakly increase
along the rows, and all entries in column $j$
belong to the set $\{1,2,\dots,\la'_j-j\}$ for
$j=1,\dots,d$. Then we define the corresponding dual flagged
Schur function $\psi_{\la/\mu}(a)$ by the formula
\ben
\psi_{\la/\mu}(a)=\sum_T \prod_{\al\in\la/\mu} a^{}_{T(\al)},
\een
summed over the dual hook $\la/\mu$-tableaux $T$.

\bth\label{thm:invtabre}
We have the expansion of the Schur function
$s_{\mu}(x)$
\ben
s_{\mu}(x)=\sum_{\la}(-1)^{m(\la/\mu)}\ts
\psi_{\la/\mu}(a)\ts \wh s_{\la}(x\vt a),
\een
summed over diagrams $\la$ which
contain $\mu$ and such that $\la$ has
$d$ boxes on the main diagonal, where
$m(\la/\mu)$ denotes the total number of boxes in
the diagram $\la/\mu$ in rows $1,\dots,d$.
\eth

\bpf
This is deduced from
Proposition~\ref{prop:inve} and the formulas for
the flagged Schur functions in \cite[8th~Variation]{m:sf}
exactly as in the proof of Theorem~\ref{thm:tabre}.
\epf

\bco\label{cor:inext}
For the expansion of the hook Schur function we have
\ben
s_{(\al\tss|\tss\be)}(x)=
\sum_{p,\ts q\geqslant 0}(-1)^p\ts
e_p(a^{}_0,a^{}_{-1},\dots,a^{}_{-\al-p+1})\ts
e_q(a^{}_1,a^{}_2,\dots,a^{}_{\be+q})\ts
\wh s^{}_{(\al+p\tss|\tss\be+q)}(x\vt a).
\een
\eco

\bex\label{ex:sonee}
We have
\ben
s_{(1)}(x)=
\sum_{p,\ts q\geqslant 0}(-1)^p\ts
a^{}_0\ts a^{}_{-1}\dots a^{}_{-p+1}\ts
a^{}_1\ts a^{}_2\dots a^{}_{q}\ts
\wh s^{}_{(p\tss|\tss q)}(x\vt a).
\een
\eex

As with the flagged Schur functions $\vp_{\la/\mu}(a)$,
we have the following factorization formula
\beql{invwhsufa}
\psi_{\la/\mu}(a)=(-1)^{m(\la/\mu)}\ts
\psi_{(\la/\mu)_-}(a)\ts\psi_{(\la'/\mu')_-}(a'),
\eeq
where $(\la/\mu)_-$ denotes the part of the skew diagram
$\la/\mu$ whose boxes lie in the rows $d+1,d+2,\dots$.
An alternative tableau presentation
for the polynomials $\psi_{(\la/\mu)_-}(a)$
is implied by the formulas \cite[$(8.2)$]{m:sf},
\cite{w:fs}. By \eqref{invwhsufa}, this also gives
alternative formulas for the
coefficients in the expansion of $s_{\mu}(x)$.

\bco\label{cor:invanexpa}
We have the tableau presentation
\ben
\psi_{(\la/\mu)_-}(a)=
\sum_T \prod_{\al\in(\la/\mu)_-} a^{}_{T(\al)},
\een
where the sum is taken
over the $(\la/\mu)_-$-tableaux $T$
whose entries in row $i$ belong to the set
$\{1,2,\dots,i-\la_i\}$ for $i=d+1,d+2,\dots$,
and the entries weakly increase along the rows
and strictly increase down
the columns.
\eco

Completing this section we note that
the canonical comultiplication on the ring $\La$
is naturally extended to the comultiplication
\ben
\Delta:\wh\La(x\vt a)\to \wh\La(x\vt a)\ot^{}_{\ts\QQ[a]}
\wh\La(x\vt a)
\een
defined on the generators by
\ben
\Delta\big(p_k(x)\big)= p_k(x)\ot 1+1\ot p_k(x).
\een
Hence, Proposition~\ref{prop:skelrp} can be interpreted
in terms of $\Delta$ as the following decomposition
of the image of the dual Schur function
\ben
\Delta\big(\wh s_{\nu}(x\vt a)\big)=
\sum_{\mu}
\wh s_{\nu/\mu}(x\vt a)
\ot \wh s_{\mu}(x\vt a)=
\sum_{\la,\ts\mu}
c_{\la\mu}^{\ts\tss\nu}(a)\ts \wh s_{\la}(x\vt a)
\ot \wh s_{\mu}(x\vt a).
\een

\section{Dual Littlewood--Richardson
polynomials}
\label{sec:dlr}
\setcounter{equation}{0}

It was pointed out in \cite[Remark~3.3]{orv:fs} that
the ring of supersymmetric functions $\La(x/y\vt a)$
is equipped with the comultiplication $\Delta$
such that
\ben
\Delta\big(p_k(x/y)\big)= p_k(x/y)\ot 1+1\ot p_k(x/y);
\een
cf. \cite[Chapter~I]{m:sfh}. The isomorphism \eqref{isom}
allows us to transfer the comultiplication
to the ring of double symmetric functions $\La(x\vt a)$
so that $\Delta$ is a $\QQ[a]$-linear ring homomorphism
\ben
\Delta:\La(x\vt a)\to \La(x\vt a)\ot^{}_{\ts\QQ[a]} \La(x\vt a)
\een
such that
\ben
\Delta\big(p_k(x\vt a)\big)= p_k(x\vt a)\ot 1+1\ot p_k(x\vt a).
\een

\bde\label{def:dulrpol}
The {\it dual Littlewood--Richardson
polynomials\/}
$\wh c_{\la\mu}^{\ts\tss\nu}(a)$ are
defined as the coefficients in the expansion
\ben
\Delta\big(s_{\nu}(x\vt a)\big)=\sum_{\la,\ts\mu}
\wh c_{\la\mu}^{\ts\tss\nu}(a)\ts s_{\la}(x\vt a)
\ot s_{\mu}(x\vt a).
\een
Equivalently, these polynomials can be found
from the decomposition
\beql{skewdec}
s_{\tss\nu/\mu}(x\vt a)=\sum_{\la}
\wh c_{\la\mu}^{\ts\tss\nu}(a)\ts s_{\la}(x\vt a).
\eeq
\vskip-1.3\baselineskip
\qed
\ede

In order to verify the equivalence of the definitions,
note that by \cite[Remark~3.3]{orv:fs},
\ben
\Delta\big(s_{\nu}(x/y\vt a)\big)=\sum_{\mu}
s_{\nu/\mu}(x/y\vt a)
\ot s_{\mu}(x/y\vt a).
\een
The claim now follows by the application
of Proposition~\ref{prop:imisom}.

It is clear from the definition that the
polynomial $\wh c_{\la\mu}^{\ts\tss\nu}(a)$ is nonzero
only if the inequality $|\nu|\geqslant |\la|+|\mu|$ holds.
In this case
it is a homogeneous polynomial in the variables $a_i$
of degree $|\nu|-|\la|-|\mu|$. Moreover, in the
particular case $|\nu|=|\la|+|\mu|$ the constant
$\wh c_{\la\mu}^{\ts\tss\nu}(a)$ equals $c_{\la\mu}^{\tss\nu}$,
the Littlewood--Richardson coefficient.

\bco\label{cor:symlr}
We have the following symmetry property
\ben
\wh c_{\la\mu}^{\ts\tss\nu}(a)=
\wh c_{\la'\mu'}^{\ts\tss\nu^{\tss\prime}}(a').
\een
\eco

\bpf
By Proposition~\ref{prop:imisom} and Definition~\ref{def:skdsf},
we have
\ben
s_{\tss\nu/\mu}(x/y\vt a)=\sum_{\la}
\wh c_{\la\mu}^{\ts\tss\nu}(a)\ts s_{\la}(x/y\vt a).
\een
The desired relations now follow from the symmetry
property \eqref{symprsu}.
\epf

We can now prove that the dual Littlewood--Richardson polynomials
$\wh c_{\la\mu}^{\ts\tss\nu}(a)$
introduced in Definition~\ref{def:dulrpol} describe
the multiplication rule for the dual Schur functions.

\bth\label{thm:proddsf}
We have
the expansion
\ben
\wh s_{\la}(x\vt a)\ts \wh s_{\mu}(x\vt a)=
\sum_{\nu} \wh c_{\la\mu}^{\ts\tss\nu}(a)\ts \wh s_{\tss\nu}(x\vt a).
\een
\eth

\bpf
We argue as in the proof of the
classical analogue of this result; see \cite[Chapter~I]{m:sfh}.
Applying Corollary~\ref{cor:ssduco} for the
families of variables
$x=x'\cup x''$ and $y=y^{\tss\prime}\cup y^{\tss\prime\prime}$ we get
\ben
\bal
\sum_{\nu\in\Pc} s_{\nu}(x/y\vt a)\ts \wh s_{\nu}(z\vt a)
{}&=\prod_{i,\ts j\geqslant 1}
\frac{1+y^{\tss\prime}_i\ts z_j}{1-x'_i\ts z_j}
\prod_{i,\ts j\geqslant 1}
\frac{1+y^{\tss\prime\prime}_i\ts z_j}{1-x''_i\ts z_j}\\
{}&=\sum_{\la,\ts\mu\in\Pc} s_{\la}(x'/y^{\tss\prime}\vt a)
\ts \wh s_{\la}(z\vt a)
\ts s_{\mu}(x''/y^{\tss\prime\prime}\vt a)\ts \wh s_{\mu}(z\vt a).
\eal
\een
On the other hand, an alternative expansion
of the sum on the left hand side is obtained
by using the relation
\ben
s_{\nu}(x/y\vt a)=\sum_{\la\subseteq\nu}s_{\la}(x'/y^{\tss\prime}\vt a)
\ts s_{\nu/\la}(x''/y^{\tss\prime\prime}\vt a)
=\sum_{\la,\ts\mu}s_{\la}(x'/y^{\tss\prime}\vt a)
\ts \wh c_{\la\mu}^{\ts\tss\nu}(a)\ts
s_{\mu}(x''/y^{\tss\prime\prime}\vt a),
\een
implied by the combinatorial formula \eqref{taba}.
Therefore, the required relation follows by
comparing the two expansions.
\epf

An explicit formula for the polynomials
$\wh c_{\la\mu}^{\ts\tss\nu}(a)$ is provided by the following
corollary, where the $c_{\al\be}^{\tss\ga}$ denote the classical
Littlewood--Richardson coefficients defined by the decomposition
of the product of the Schur functions
\ben
s_{\al}(x)\ts s_{\be}(x)=
\sum_{\ga} c_{\al\be}^{\tss\ga}\ts s_{\tss\ga}(x).
\een
We suppose that $\vp_{\al/\la}(a)=\psi_{\al/\la}(a)=0$
unless $\la\subseteq\al$ and the diagrams $\la$ and $\al$ have
the same number of boxes on their main diagonals.

\bco\label{cor:dulr}
We have
\ben
\wh c_{\la\mu}^{\ts\tss\nu}(a)=\sum_{\al,\be,\ga}
(-1)^{n(\al/\la)+n(\be/\mu)+m(\nu/\ga)}\ts
c_{\al\be}^{\tss\ga}\ts
\vp_{\al/\la}(a)\ts \vp_{\be/\mu}(a)\ts
\psi_{\nu/\ga}(a),
\een
summed over diagrams $\al$, $\be$, $\ga$.
In particular,
$\wh c_{\la\mu}^{\ts\tss\nu}(a)=0$ unless $\la\subseteq\nu$
and $\mu\subseteq\nu$.
\eco

\bpf
The formula follows from Theorems~\ref{thm:tabre}, \ref{thm:invtabre}
and \ref{thm:proddsf}. The second statement is implied
by the same property of the Littlewood--Richardson coefficients.
\epf

\bex\label{ex:skewsf}
If $k\leqslant l$ and $k+l\leqslant m$ then
\ben
\wh c_{(k)(l)}^{\ts\tss(m)}(a)=\sum_{r+s=m-k-l}
(-1)^s\ts h_r(a_0,a_{-1},\dots,a_{-k+1})
\ts e_s(a_{-l},a_{-l-1},\dots,a_{-m+2}).
\een
In particular,
\ben
\wh c_{(1)(l)}^{\ts\tss(m)}(a)=
(a_0-a_{-l})(a_0-a_{-l-1})\dots(a_0-a_{-m+2}).
\een
Applying Corollary~\ref{cor:symlr}, we also get
\ben
\wh c_{(1^k)(1^l)}^{\ts\tss(1^m)}(a)=\sum_{r+s=m-k-l}
(-1)^r\ts h_r(a_1,a_2,\dots,a_k)
\ts e_s(a_{l+1},a_{l+2},\dots,a_{m-1})
\een
and
\ben
\wh c_{(1)(1^l)}^{\ts\tss(1^m)}(a)=
(a_{l+1}-a_1)(a_{l+2}-a_1)\dots(a_{m-1}-a_1).
\een
These relations provide explicit formulas for the images
of the double elementary and complete symmetric
functions $h_m(x\vt a)$ and $e_m(x\vt a)$ with respect to
the comultiplication $\Delta$.
\qed
\eex

Another formula for the
dual Littlewood--Richardson polynomials
$\wh c_{\la\mu}^{\ts\tss\nu}(a)$ can be obtained with the use
of the decomposition \eqref{skewdec}. We will consider the
skew double Schur function as the sequence of polynomials
$s_{\nu/\mu}(x\vt a)$ defined in \eqref{seqpol}.
For a given skew diagram $\nu/\mu$ consider the finite set
of variables $x=(x_1,\dots,x_n)$, where
$\nu'_j-\mu'_j\leqslant n$ for all $j$; that is,
the number of boxes in each column of $\nu/\mu$ does not
exceed $n$. Since the skew double Schur functions are
consistent with the evaluation homomorphisms \eqref{eval},
the polynomials $\wh c_{\la\mu}^{\ts\tss\nu}(a)$
are determined by the decomposition \eqref{skewdec},
where $x$ is understood as the above finite set of variables.

In order to formulate the result, introduce
$\nu/\mu$-{\it supertableaux\/} $T$
which are obtained by filling in the boxes of
$\nu/\mu$ with the symbols $1,1^{\tss\prime},\dots,n,n'$
in such a way that
in each row (resp. column) each primed
index is to the left (resp. above)
of each unprimed index;
unprimed indices weakly decrease along the rows and
strictly decrease down the columns;
primed indices strictly increase along the rows and
weakly increase down the columns.

Introduce the ordering on the set of boxes of
a skew diagram by reading them by columns from
left to right and from bottom to top in each column.
We call this the {\it column order\/}.
We shall write $\al\prec \be$ if $\al$ (strictly) precedes $\be$
with respect to the column order.

Suppose that $\la$ is a diagram. Given a sequence
of diagrams $R$ of the form
\beql{r}
\varnothing=\rho^{(0)}\to\rho^{(1)}\to
\dots\to\rho^{(l-1)}\to\rho^{(l)}=\la,
\eeq
we let $r_i$ denote the row number of the box added to
the diagram $\rho^{(i-1)}$.
The sequence $r_1r_2\dots r_l$ is called the {\it Yamanouchi symbol\/}
of $R$.  Construct the set $\Tc(\nu/\mu,R)$ of {\it barred\/}
$\nu/\mu$-supertableaux
$T$ such that $T$ contains
boxes $\al_1,\dots,\al_l$ with
\ben
\al_1\prec\dots\prec\al_l\Fand T(\al_i)=r_i,\quad
1\leqslant i\leqslant l,
\een
where all entries $r_i$ are unprimed and the boxes
are listed in the column order
which is restricted to the
subtableau of $T$ formed by the unprimed indices.

We will distinguish the entries in $\al_1,\dots,\al_l$ by barring
each of them. So, an element of $\Tc(\nu/\mu,R)$ is a pair consisting
of a $\nu/\mu$-supertableau and a chosen sequence of barred entries
compatible with $R$.
We shall keep the notation $T$ for such a pair.

For each box $\alpha$ with $\al_i\prec\al\prec\al_{i+1}$,
$0\leqslant i\leqslant l$,
which is occupied by an unprimed index,
set $\rho(\alpha)=\rho^{(i)}$.

\bth\label{thm:supt}
The dual Littlewood--Richardson polynomials
can be given by
\ben
\bal
\wh c_{\la\mu}^{\ts\tss\nu}(a)=\sum_R
\sum_{T}
\prod_{\underset{\scriptstyle T(\alpha)
\text{\ts\ts unprimed,\ts unbarred}}{\alpha\in\nu/\mu}}
&\big(a^{}_{\tss T(\alpha)-\rho(\alpha)^{}_{T(\al)}}-
a^{}_{\tss T(\alpha)-c(\alpha)}\big)\\
{}\times{}\prod_{\underset{\scriptstyle T(\alpha)
\text{\ts\ts primed}}{\alpha\in\nu/\mu}}
&\big(a^{}_{\tss T(\alpha)-c(\alpha)}
-a^{}_{\tss T(\alpha)}\big),
\eal
\een
summed over sequences $R$ of the form \eqref{r} and
barred supertableaux $T\in\Tc(\nu/\mu,R)$.
\eth

\bpf
Due to \eqref{seqpol}, we have
\ben
\wh c_{\la\mu}^{\ts\tss\nu}(a)
=\sum_{\mu\subseteq\ts\rho\ts\subseteq\tss\nu}
\wt c_{\la\rho}^{\ts\tss\nu}(a)\ts
s_{\rho^{\tss\prime}/\mu^{\tss\prime}}(-a^{(n)} \tss|\tss {-}a),
\een
where the polynomials $\wt c_{\la\rho}^{\ts\tss\nu}(a)$
are defined by the decomposition
\ben
\wt s_{\nu/\rho}(x\vt a)=\sum_{\la}
\wt c_{\la\rho}^{\ts\tss\nu}(a)\ts s_{\la}(x\vt a).
\een
The desired formula is now implied by
\cite[Lemma~2.4]{m:lr} which gives the combinatorial
expression for the coefficients
$\wt c_{\la\rho}^{\ts\tss\nu}(a)$ and thus takes care
of the unprimed part of $T$; the expression for the primed part
is implied by \eqref{facsf}.
\epf

\bre\label{rem:open}
Both the formulas for $\wh c_{\la\mu}^{\ts\tss\nu}(a)$
provided by Corollary~\ref{cor:dulr}
and Theorem~\ref{thm:supt}
involve some terms which cancel pairwise.
It would be interesting to find a combinatorial presentation
of the polynomials $\wh c_{\la\mu}^{\ts\tss\nu}(a)$
analogous to \cite{kt:pe}, \cite{k:elr} or \cite{m:lr}
and to understand their positivity properties.
A possible way to find such a presentation could rely on
the vanishing theorem of the supersymmetric Schur functions
obtained in \cite[Theorems~5.1 \& 5.2]{orv:fs}; see also
\cite[Theorem~4.4]{m:fs} for a similar result.
\qed
\ere

\bex\label{ex:superca}
In order to calculate the polynomial
$\wh c_{(1)\ts(2)}^{\ts\tss(2^2)}(a)$, take
$\la=(1)$, $\mu=(2)$, $\nu=(2^2)$ and
$n=1$. The barred
supertableaux compatible with the sequence $\varnothing\to (1)$
are

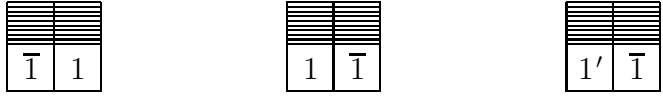
\setlength{\unitlength}{0.75em}
\begin{center}
\begin{picture}(30,4.2)

\put(0,0){\line(0,1){3.8}}
\put(2,0){\line(0,1){3.8}}
\put(4,0){\line(0,1){3.8}}
\put(0,0){\line(1,0){4}}
\put(0,2){\line(1,0){4}}
\put(0.7,0.5){$\overline{1}$}
\put(0.7,2.5){}
\put(2.7,2.5){}
\put(2.7,0.5){1}
\put(0,2.2){\line(1,0){4}}
\put(0,2.4){\line(1,0){4}}
\put(0,2.6){\line(1,0){4}}
\put(0,2.8){\line(1,0){4}}
\put(0,3.0){\line(1,0){4}}
\put(0,3.2){\line(1,0){4}}
\put(0,3.4){\line(1,0){4}}
\put(0,3.6){\line(1,0){4}}
\put(0,3.8){\line(1,0){4}}

\put(12,0){\line(0,1){3.8}}
\put(14,0){\line(0,1){3.8}}
\put(16,0){\line(0,1){3.8}}
\put(12,0){\line(1,0){4}}
\put(12,2){\line(1,0){4}}
\put(12.7,0.5){1}
\put(12.7,2.5){}
\put(14.7,2.5){}
\put(14.7,0.5){$\overline{1}$}
\put(12,2.2){\line(1,0){4}}
\put(12,2.4){\line(1,0){4}}
\put(12,2.6){\line(1,0){4}}
\put(12,2.8){\line(1,0){4}}
\put(12,3.0){\line(1,0){4}}
\put(12,3.2){\line(1,0){4}}
\put(12,3.4){\line(1,0){4}}
\put(12,3.6){\line(1,0){4}}
\put(12,3.8){\line(1,0){4}}

\put(24,0){\line(0,1){3.8}}
\put(26,0){\line(0,1){3.8}}
\put(28,0){\line(0,1){3.8}}
\put(24,0){\line(1,0){4}}
\put(24,2){\line(1,0){4}}
\put(24.5,0.5){$1^{\tss\prime}$}
\put(24.7,2.5){}
\put(26.7,2.5){}
\put(26.7,0.5){$\overline{1}$}
\put(24,2.2){\line(1,0){4}}
\put(24,2.4){\line(1,0){4}}
\put(24,2.6){\line(1,0){4}}
\put(24,2.8){\line(1,0){4}}
\put(24,3.0){\line(1,0){4}}
\put(24,3.2){\line(1,0){4}}
\put(24,3.4){\line(1,0){4}}
\put(24,3.6){\line(1,0){4}}
\put(24,3.8){\line(1,0){4}}

\end{picture}
\end{center}
\setlength{\unitlength}{1pt}

\noindent
so that
\ben
\wh c_{(1)\ts(2)}^{\ts\tss(2^2)}(a)
=a_0-a_1+a_1-a_2
+a_{2}-a_1=a_0-a_1.
\een
Alternatively, we can take
$\la=(2)$, $\mu=(1)$, $\nu=(2^2)$ and
$n=2$. The barred
supertableaux compatible with the sequence $\varnothing\to (1)\to (2)$
are

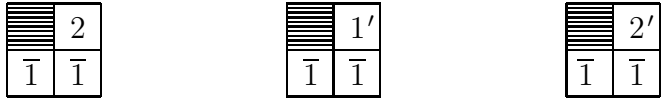
\setlength{\unitlength}{0.75em}
\begin{center}
\begin{picture}(30,4.6)

\put(0,0){\line(0,1){4}}
\put(2,0){\line(0,1){4}}
\put(4,0){\line(0,1){4}}
\put(0,0){\line(1,0){4}}
\put(0,2){\line(1,0){4}}
\put(0,4){\line(1,0){4}}
\put(0.7,0.5){$\overline{1}$}
\put(0.7,2.5){}
\put(2.7,2.5){2}
\put(2.7,0.5){$\overline{1}$}
\put(0,2.2){\line(1,0){2}}
\put(0,2.4){\line(1,0){2}}
\put(0,2.6){\line(1,0){2}}
\put(0,2.8){\line(1,0){2}}
\put(0,3.0){\line(1,0){2}}
\put(0,3.2){\line(1,0){2}}
\put(0,3.4){\line(1,0){2}}
\put(0,3.6){\line(1,0){2}}
\put(0,3.8){\line(1,0){2}}

\put(12,0){\line(0,1){4}}
\put(14,0){\line(0,1){4}}
\put(16,0){\line(0,1){4}}
\put(12,0){\line(1,0){4}}
\put(12,2){\line(1,0){4}}
\put(12,4){\line(1,0){4}}
\put(12.7,0.5){$\overline{1}$}
\put(12.7,2.5){}
\put(14.7,2.5){$1^{\tss\prime}$}
\put(14.7,0.5){$\overline{1}$}
\put(12,2.2){\line(1,0){2}}
\put(12,2.4){\line(1,0){2}}
\put(12,2.6){\line(1,0){2}}
\put(12,2.8){\line(1,0){2}}
\put(12,3.0){\line(1,0){2}}
\put(12,3.2){\line(1,0){2}}
\put(12,3.4){\line(1,0){2}}
\put(12,3.6){\line(1,0){2}}
\put(12,3.8){\line(1,0){2}}

\put(24,0){\line(0,1){4}}
\put(26,0){\line(0,1){4}}
\put(28,0){\line(0,1){4}}
\put(24,0){\line(1,0){4}}
\put(24,2){\line(1,0){4}}
\put(24,4){\line(1,0){4}}
\put(24.5,0.5){$\overline{1}$}
\put(24.7,2.5){}
\put(26.7,2.5){$2^{\tss\prime}$}
\put(26.7,0.5){$\overline{1}$}
\put(24,2.2){\line(1,0){2}}
\put(24,2.4){\line(1,0){2}}
\put(24,2.6){\line(1,0){2}}
\put(24,2.8){\line(1,0){2}}
\put(24,3.0){\line(1,0){2}}
\put(24,3.2){\line(1,0){2}}
\put(24,3.4){\line(1,0){2}}
\put(24,3.6){\line(1,0){2}}
\put(24,3.8){\line(1,0){2}}

\end{picture}
\end{center}
\setlength{\unitlength}{1pt}

\noindent
so that
\ben
\wh c_{(1)\ts(2)}^{\ts\tss(2^2)}(a)
=a_2-a_1+a_0-a_1
+a_{1}-a_2=a_0-a_1.
\een
This agrees with the previous calculation
and the formula implied by Corollary~\ref{cor:dulr}.
\eex

\bex\label{ex:supediff}
Theorem~\ref{thm:supt} gives formulas
for the polynomials $\wh c_{(k)(l)}^{\ts\tss(m)}(a)$
and $\wh c_{(1^k)(1^l)}^{\ts\tss(1^m)}(a)$
in a different form
as compared to Example~\ref{ex:skewsf}.
If $k+l\leqslant m$ then
\ben
\bal
\wh c_{(k)(l)}^{\ts\tss(m)}(a)=\sum
&(a_0-a_{-l})(a_0-a_{-l-1})\cdots (a_0-a_{-l-i_1+1})\\
&{}\times{}(a_{-1}-a_{-l-i_1-1})\cdots (a_{-1}-a_{-l-i_2+1})\\
&{}\times{}\cdots
(a_{-k+1}-a_{-l-i_{k-1}-1})\cdots (a_{-k+1}-a_{-m+2})
\eal
\een
summed over the sets of indices
$0\leqslant i_1<\dots<i_{k-1}\leqslant m-l-2$.
A similar expression for $\wh c_{(1^k)(1^l)}^{\ts\tss(1^m)}(a)$
follows by the application of Corollary~\ref{cor:symlr}.
\eex

\section{Transition matrices}
\label{sec:ome}
\setcounter{equation}{0}

\subsection{Pairing between the double and dual symmetric
functions}

We now prove alternative expansion formulas for the infinite product
which occurs in the Cauchy formula \eqref{cchid}.
These formulas turn into the well known identities
when $a$ is specialized to the sequence of zeros; see
\cite[Chapter~I]{m:sfh}.

Let $\la=(\la_1,\dots,\la_l)$ be a partition and suppose
that the length of $\la$ does not
exceed $l$. Using the notation
\eqref{dumon}, introduce the
{\it dual monomial symmetric function\/}
$\wh m_{\la}(x\vt a)\in\wh\La(x\vt a)$ by the formula
\ben
\wh m_{\la}(x\vt a)=
\sum_{\si} (x_{\si(1)},a)^{\la_1}\ts (x_{\si(2)},a)^{\la_2}\dots
(x_{\si(l)},a)^{\la_l},
\een
summed over permutations $\si$ of the $x_i$
which give distinct monomials.

For a partition $\la=(1^{m_1}\ts 2^{m_2}\dots)$ set
$z_{\la}=\prod_{i\geqslant 1} i^{m_i}\ts m_i!$.

\bpr\label{prop:expa}
We have
the expansions
\beql{hmde}
\prod_{i,\ts j\geqslant 1}\frac{1-a_i\ts y_j}{1-x_i\ts y_j}
=\sum_{\la\in\Pc}
h_{\la}(x\vt a)\ts \wh m_{\la}(y\vt a)
\eeq
and
\beql{ppde}
\prod_{i,\ts j\geqslant 1}\frac{1-a_i\ts y_j}{1-x_i\ts y_j}
=\sum_{\la\in\Pc}
z_{\la}^{-1}\ts p_{\la}(x\vt a)\ts p_{\la}(y).
\eeq
\epr

\bpf
Let us set
\ben
H(t)=\prod_{i=1}^{\infty}\frac{1-a_i\tss t}{1-x_i\tss t}.
\een
Then using \eqref{genh} and arguing as in
\cite[Chapter~I]{m:sfh}, we can write
\ben
\prod_{i,\ts j\geqslant 1}\frac{1-a_i\ts y_j}{1-x_i\ts y_j}
=\prod_{j\geqslant 1} H(y_j)
=\prod_{j\geqslant 1}\sum_{k=0}^{\infty} h_k(x\vt a)
\ts (y_j,a)^k=\sum_{\la\in\Pc}
h_{\la}(x\vt a)\ts \wh m_{\la}(y\vt a),
\een
which proves \eqref{hmde}. For the proof of \eqref{ppde}
note that
\ben
\bal
\ln H(t)&=\sum_{i\geqslant 1}\Big(\ln(1-a_it)-\ln(1-x_it)\Big)\\
{}&=\sum_{i\geqslant 1}\sum_{k\geqslant 1}
\Big(\frac{x_i^k\ts t^k}{k}-\frac{a_i^k\ts t^k}{k}\Big)
=\sum_{k\geqslant 1}\frac{p_k(x\vt a) t^k}{k}.
\eal
\een
Hence,
\ben
H(t)=\sum_{\la\in\Pc} z_{\la}^{-1}\ts p_{\la}(x\vt a)\ts t^{|\la|}.
\een
Now apply this relation to the sets of variables
$x$ and $a$ respectively replaced with the sets $\{x_iy_j\}$
and $\{a_iy_j\}$. Then $p_k(x\vt a)$ is replaced by
$p_{\la}(x\vt a)\ts p_{\la}(y)$, and \eqref{ppde}
follows by putting $t=1$.
\epf

Now define the $\QQ[a]$-bilinear pairing
between the rings $\La(x\vt a)$ and $\wh\La(y\vt a)$,
\beql{pairi}
\langle\ ,\ \rangle:\big(\La(x\vt a),\wh\La(y\vt a)\big)\to \QQ[a],
\eeq
by setting
\beql{defpa}
\big\langle h_{\la}(x\vt a),\wh m_{\mu}(y\vt a)\big\rangle
=\de_{\la\mu}.
\eeq
Clearly, $\langle u,\wh v\rangle$ is
a well-defined polynomial in $a$ for any elements
$u\in\La(x\vt a)$ and $\wh v\in\wh\La(y\vt a)$ which is determined
from \eqref{defpa} by linearity.

The following is an analogue of the duality properties
of the classical bases of the ring of symmetric functions;
see \cite[Chapter~I]{m:sfh}.

\bpr\label{prop:dubas}
Let $\{u_{\la}(x\vt a)\}$
and $\{\wh v_{\la}(y\vt a)\}$ be families of elements
of rings $\La(x\vt a)$ and $\wh\La(y\vt a)$, respectively,
which are parameterized by
all partitions. Suppose that for any $n\geqslant 0$
the highest degree components in $x$ {\rm(}resp.,
the lowest degree components in $y${\rm)}
of the elements $u_{\la}(x\vt a)$ {\rm(}resp.,
$\wh v_{\la}(y\vt a)${\rm)}
with $|\la|=n$
form a basis of the space of homogeneous symmetric
functions in $x$ {\rm(}resp., $y${\rm)} of degree $n$.
Then the following conditions are equivalent:
\beql{uvdu}
\big\langle u_{\la}(x\vt a),\wh v_{\mu}(y\vt a)\big\rangle
=\de_{\la\mu},\qquad\text{for all}\quad \la,\mu;
\eeq
\beql{unsum}
\sum_{\la\in\Pc}
u_{\la}(x\vt a)\ts \wh v_{\la}(y\vt a)
=\prod_{i,\ts j\geqslant 1}\frac{1-a_i\ts y_j}{1-x_i\ts y_j}.
\eeq
\epr

\bpf
We only need to slightly modify the respective argument of
\cite[Chapter~I]{m:sfh}.
Write
\ben
u_{\la}(x\vt a)=\sum_{\rho} A_{\la\rho}(a)\ts h_{\rho}(x\vt a),
\qquad
\wh v_{\mu}(y\vt a)=\sum_{\si} B_{\mu\si}(a)\ts \wh m_{\si}(y\vt a),
\een
where the first sum is taken over partitions $\rho$
with $|\rho|\leqslant |\la|$, while the second
is taken over partitions $\si$
with $|\si|\geqslant |\mu|$.
Then
\ben
\big\langle u_{\la}(x\vt a),\wh v_{\mu}(y\vt a)\big\rangle
=\sum_{\rho} A_{\la\rho}(a)\ts B_{\mu\rho}(a).
\een
Hence, condition \eqref{uvdu} is equivalent to
\beql{abmu}
\sum_{\rho} A_{\la\rho}(a)\ts B_{\mu\rho}(a)=\de_{\la\mu}.
\eeq
On the other hand, due to \eqref{hmde},
\eqref{unsum} can be written as
\ben
\sum_{\la\in\Pc}
u_{\la}(x\vt a)\ts \wh v_{\la}(y\vt a)=\sum_{\rho\in\Pc}
h_{\rho}(x\vt a)\ts \wh m_{\rho}(y\vt a),
\een
which is equivalent to
\ben
\sum_{\la} A_{\la\rho}(a)\ts B_{\la\si}(a)=\de_{\rho\si}.
\een
This condition is easily verified
to be equivalent to \eqref{abmu}.
\epf

Applying Theorem~\ref{thm:cch} and Proposition~\ref{prop:expa} we get
the following corollary.

\bco\label{cor:dumup}
Under the pairing \eqref{pairi} we have
\ben
\big\langle s_{\la}(x\vt a),\wh s_{\mu}(y\vt a)\big\rangle
=\de_{\la\mu}
\Fand
\big\langle p_{\la}(x\vt a),p_{\mu}(y)\big\rangle
=\de_{\la\mu}\ts z_{\la}.
\een
\vskip-1.2\baselineskip
\qed
\eco

Thus, the symmetric functions $\wh s_{\la}(y\vt a)$
are dual to the double Schur functions
$s_{\la}(x\vt a)$
in sense of the pairing \eqref{pairi}.

Using the isomorphism \eqref{isom} and the pairing
\eqref{pairi}, we get another
$\QQ[a]$-bilinear pairing
\beql{pairisu}
\langle\ ,\ \rangle:\big(\La(x/y\vt a),\wh\La(z\vt a)\big)\to \QQ[a]
\eeq
such that
\beql{duoda}
\big\langle s_{\la}(x/y\vt a),\wh s_{\mu}(z\vt a)\big\rangle
=\de_{\la\mu}.
\eeq
Note that Proposition~\ref{prop:dubas} can be easily
reformulated for the pairing \eqref{pairisu}. In particular,
the condition \eqref{unsum} is now replaced by
\beql{unsumsu}
\sum_{\la\in\Pc}
u_{\la}(x/y\vt a)\ts \wh v_{\la}(z\vt a)
=\prod_{i,\ts j\geqslant 1}\frac{1+y_i\ts z_j}{1-x_i\ts z_j}.
\eeq
This implies that
\beql{duod}
\big\langle s_{\la}(x/y),s_{\mu}(z)\big\rangle
=\de_{\la\mu},
\eeq
where $s_{\la}(x/y)$ denotes the ordinary supersymmetric
Schur function which is obtained from $s_{\la}(x/y\vt a)$
by the specialization $a_i=0$.
Together with Theorems~\ref{thm:tabre} and
\ref{thm:invtabre}, the relations
\eqref{duoda} and \eqref{duod} imply the following
expansions for the supersymmetric Schur functions.

\bco\label{cor:expasdo}
We have the decompositions
\ben
s_{\la}(x/y\vt a)=\sum_{\mu}(-1)^{m(\la/\mu)}\ts
\psi_{\la/\mu}(a)\ts s_{\mu}(x/y),
\een
summed over diagrams $\mu$ contained
in $\la$ and such that $\la$ and $\mu$ have
the same number of
boxes on the main diagonal; and
\ben
s_{\mu}(x/y)=\sum_{\la}(-1)^{n(\la/\mu)}\ts
\vp_{\la/\mu}(a)\ts s_{\la}(x/y\vt a),
\een
summed over diagrams $\la$ which
contain $\mu$ and such that $\la$ and $\mu$ have
the same number of
boxes on the main diagonal.
\qed
\eco

Note that expressions for $\psi_{\la/\mu}(a)$ and
$\vp_{\la/\mu}(a)$ in terms of determinants as
in Propositions~\ref{prop:dsfexp} and \ref{prop:inve}
were given in \cite{orv:fs}.
Corollary~\ref{cor:expasdo}
gives new tableau formulas for these coefficients.
Moreover, under the specialization $a_i=-i+1/2$
the supersymmetric Schur functions
$s_{\la}(x/y\vt a)$ turn into the {\it Frobenius--Schur
functions\/} $Fs_{\mu}$; see \cite{orv:fs}.
Hence, the transition coefficients between the
$Fs_{\mu}$ and the Schur functions can be
found as follows; cf. \cite[Theorem~2.6]{orv:fs}.

\bco\label{cor:frschur}
We have the decompositions
\ben
Fs_{\la}=\sum_{\mu}(-1)^{m(\la/\mu)}\ts
\psi_{\la/\mu}\ts s_{\mu}(x/y)
\een
and
\ben
s_{\mu}(x/y)=\sum_{\la}(-1)^{n(\la/\mu)}\ts
\vp_{\la/\mu}\ts Fs_{\la},
\een
where $\psi_{\la/\mu}$
and $\vp_{\la/\mu}$ are the respective values of
the polynomials $\psi_{\la/\mu}(a)$ and
$\vp_{\la/\mu}(a)$ at $a_i=-i+1/2$,
$i\in\ZZ$.
\qed
\eco

Using the notation of Corollary~\ref{cor:expasdo} and
applying the isomorphism \eqref{isom} we get the respective
expansion formulas involving the double Schur functions.

\bco\label{cor:expasdouble}
We have the decompositions
\ben
s_{\la}(x\vt a)=\sum_{\mu}(-1)^{m(\la/\mu)}\ts
\psi_{\la/\mu}(a)\ts s_{\mu}(x/{-}a^+)
\een
and
\ben
s_{\mu}(x/{-}a^+)=\sum_{\la}(-1)^{n(\la/\mu)}\ts
\vp_{\la/\mu}(a)\ts s_{\la}(x\vt a),
\een
where $a^+=(a_1,a_2,\dots)$.
\qed
\eco

Expressions for the coefficients
in the expansions relating the double and ordinary
Schur functions or polynomials can be found
in \cite{k:elr}, \cite{k:pf}, \cite{m:sf},
\cite{m:lr} and \cite{ms:lr}.

Let us now recall
the isomorphism $\om_a:\La(x\vt a)\to\La(x\vt a')$ and
the involution $\wh\om:\wh\La(x\vt a)\to \wh\La(x\vt a)$;
see \eqref{omega} and \eqref{whome}.
Since every polynomial $c(a)\in \QQ[a]$ can be regarded
as an element of $\QQ[a']$, the ring $\wh\La(x\vt a')$
can be naturally identified with $\wh\La(x\vt a)$
via the map $c(a)\mapsto c^{\tss\prime}(a')$, where
$c^{\tss\prime}(a')=c(a)$
as polynomials in the $a_i$, $i\in\ZZ$.

\bpr\label{prop:oma}
For any elements $u\in \La(x\vt a)$ and $\wh v\in \wh\La(y\vt a)$
we have
\ben
\big\langle \om_a\tss u,\ts\wh\om\tss \wh v\big\rangle'
=\big\langle u, \wh v\big\rangle,
\een
where $\langle \ ,\ \rangle'$ denotes the pairing
\eqref{pairi} between $\La(x\vt a')$
and $\wh\La(y\vt a)\simeq \wh\La(y\vt a')$.
\epr

\bpf
It suffices to take $u=s_{\la}(x\vt a)$ and
$\wh v=\wh s_{\mu}(y\vt a)$.
Using \eqref{omsla} and \eqref{whome}, we get
\ben
\big\langle \om_a\tss s_{\la}(x\vt a),
\ts\wh\om\tss \wh s_{\mu}(y\vt a)\big\rangle'
=\big\langle s_{\la'}(x\vt a'),
\ts \wh s_{\mu'}(y\vt a')\big\rangle'.
\een
By Corollary~\ref{cor:dumup} this equals $\de_{\la\mu}$,
and hence coincides with
$\langle s_{\la}(x\vt a),
\ts \wh s_{\mu}(y\vt a)\rangle$.
\epf

Introduce the {\it dual forgotten symmetric functions\/}
$\wh f_{\la}(y\vt a)\in\wh\La(y\vt a)$ as the images
of the dual monomial symmetric functions
under the involution $\wh\om$, that is,
\ben
\wh f_{\la}(y\vt a)=\wh\om\ts \wh m_{\la}(y\vt a'),\qquad \la\in\Pc.
\een
Furthermore,
for any partition $\la$ define the {\it double
monomial symmetric functions\/} $m_{\la}(x\vt a)\in\La(x\vt a)$
and the {\it double forgotten symmetric functions\/}
$f_{\la}(x\vt a)\in\La(x\vt a)$ by the relations
\beql{hmdua}
\prod_{i,\ts j\geqslant 1}\frac{1-a_i\ts y_j}{1-x_i\ts y_j}
=\sum_{\la\in\Pc}
m_{\la}(x\vt a)\ts \wh h_{\la}(y\vt a)
\eeq
and
\beql{ppdua}
\prod_{i,\ts j\geqslant 1}\frac{1-a_i\ts y_j}{1-x_i\ts y_j}
=\sum_{\la\in\Pc}
f_{\la}(x\vt a)\ts \wh e_{\la}(y\vt a).
\eeq

Hence, by Proposition~\ref{prop:dubas},
under the pairing \eqref{pairi}
we have
\beql{pamh}
\big\langle m_{\la}(x\vt a),\wh h_{\mu}(y\vt a)\big\rangle
=\de_{\la\mu}
\Fand
\big\langle f_{\la}(x\vt a),\wh e_{\mu}(y\vt a)\big\rangle
=\de_{\la\mu}.
\eeq
Moreover, Propositions~\ref{prop:om} and \ref{prop:oma} imply
\ben
\om_a: m_{\la}(x\vt a)\mapsto h_{\la}(x\vt a'),
\qquad f_{\la}(x\vt a)\mapsto e_{\la}(x\vt a'),
\qquad p_{\la}(x\vt a)\mapsto \ve_{\la}\ts p_{\la}(x\vt a'),
\een
where $\ve_{\la}=(-1)^{|\la|-\ell(\la)}$. To check the latter
relation we need to recall that under
the involution $\om$ of the ring of symmetric functions we have
$\om:p_{\la}(y)\mapsto
\ve_{\la}\ts p_{\la}(y)$; see \cite[Chapter~I]{m:sfh}.

We can now obtain
analogues of the decomposition of Corollary~\ref{cor:duco}
for other families of symmetric functions.

\bco\label{cor:decmf}
We have the decompositions
\ben
\bal
\prod_{i,\ts j\geqslant 1}\frac{1+x_i\ts y_j}{1+a_i\ts y_j}
&=\sum_{\la\in\Pc} e_{\la}(x\vt a)\ts
\wh m_{\la}(y\vt a^{\ts\prime}),\\[1em]
\prod_{i,\ts j\geqslant 1}\frac{1+x_i\ts y_j}{1+a_i\ts y_j}
&=\sum_{\la\in\Pc} \ve_{\la}\ts z_{\la}^{-1}\ts
p_{\la}(x\vt a)\ts p_{\la}(y),\\[1em]
\prod_{i,\ts j\geqslant 1}\frac{1+x_i\ts y_j}{1+a_i\ts y_j}
&=\sum_{\la\in\Pc} m_{\la}(x\vt a)\ts \wh e_{\la}(y\vt a^{\ts\prime}).
\eal
\een
\eco

\bpf
The relations follow by the application of $\om_a$
to the expansions \eqref{hmde}, \eqref{ppde} and
by the application of $\wh\om$ to \eqref{hmdua}.
\epf

Note that relations of this kind involving the forgotten
symmetric functions can be obtained in a similar way.

\subsection{Kostka-type and character polynomials}

The entries of the transition matrices between the classical
bases of the ring of symmetric functions
can be expressed in terms of the Kostka numbers $K_{\la\mu}$
and the values $\chi^{\la}_{\mu}$ of the irreducible
characters of the symmetric groups; see \cite[Sections~I.6,~I.7]{m:sfh}.
By analogy with the classical case, introduce the
{\it Kostka-type polynomials\/} $K_{\la\mu}(a)$ and
the {\it character polynomials\/} $\chi^{\la}_{\mu}(a)$
as well as their dual counterparts
$\wh K_{\la\mu}(a)$ and $\wh\chi^{\ts\la}_{\mu}(a)$
by the respective expansions
\ben
s_{\la}(x\vt a)=\sum_{\mu} K_{\la\mu}(a)\ts m_{\mu}(x\vt a),
\qquad
\wh s_{\la}(y\vt a)=\sum_{\mu}
\wh K_{\la\mu}(a)\ts \wh m_{\mu}(y\vt a),
\een
and
\ben
p_{\mu}(x\vt a)=\sum_{\la} \chi^{\la}_{\mu}(a)\ts s_{\la}(x\vt a),
\qquad
p_{\mu}(y)=\sum_{\la} \wh\chi^{\ts\la}_{\mu}(a)\ts \wh s_{\la}(y\vt a).
\een
If $|\la|=|\mu|$, then
\beql{cla}
K_{\la\mu}(a)=\wh K_{\la\mu}(a)=K_{\la\mu}\Fand
\chi^{\la}_{\mu}(a)=\wh\chi^{\ts\la}_{\mu}(a)=\chi^{\la}_{\mu}.
\eeq
Moreover, $K_{\la\mu}(a)$ and $\wh\chi^{\ts\la}_{\mu}(a)$
are zero unless $|\la|\geqslant |\mu|$, while
$\wh K_{\la\mu}(a)$ and $\chi^{\la}_{\mu}(a)$
are zero unless $|\la|\leqslant |\mu|$.

Using the duality properties of the double and dual symmetric
functions, we can get all other transition matrices
in the same way as this is done in \cite[Sections~I.6,~I.7]{m:sfh}.
In particular, we have the relations
\ben
h_{\mu}(x\vt a)=\sum_{\la} \wh K_{\la\mu}(a)\ts s_{\la}(x\vt a),
\qquad
\wh h_{\mu}(y\vt a)=\sum_{\la} K_{\la\mu}(a)\ts \wh s_{\la}(y\vt a).
\een

The Littlewood--Richardson
polynomials $c_{\la\mu}^{\tss\nu}(a)$ defined in \eqref{lrpoldef}
are Graham positive
as they can be written as
polynomials in the differences $a_i-a_j$, $i<j$, with
positive integer coefficients; see \cite{g:pe}.
Explicit positive formulas for
$c_{\la\mu}^{\ts\nu}(a)$
were found in \cite{kt:pe}, \cite{k:elr} and \cite{m:lr}.
Using the fact that $h_k(x\vt a)$ coincides with
$s_{(k)}(x\vt a)$, we come to the following
expression for the polynomials $\wh K_{\la\mu}(a)$:
\ben
\wh K_{\la\mu}(a)=\sum_{\rho^{(1)},\dots,\ts\rho^{(l-2)}}
c_{(\mu_1)\tss\rho^{(1)}}^{\tss\la}(a)\ts
c_{(\mu_2)\tss\rho^{(2)}}^{\tss\rho^{(1)}}(a)\ts
\dots c_{(\mu_{l-2})\tss\rho^{(l-2)}}^{\tss\rho^{(l-3)}}(a)\ts
c_{(\mu_{l-1})\tss(\mu_{l})}^{\tss\rho^{(l-2)}}(a),
\een
summed over partitions $\rho^{(i)}$,
where $\mu=(\mu_1,\dots\mu_l)$. In particular, each
dual Kostka-type polynomial
$\wh K_{\la\mu}(a)$ is Graham positive. For a more explicit
tableau presentation of the polynomials $\wh K_{\la\mu}(a)$
see \cite{f:tm}.

\bex\label{ex:kona}
We have
\ben
\wh K_{(3\ts 2)\tss(3\ts 2\ts 1)}(a)=
\sum_{\rho}
c_{(3)\tss\rho}^{\tss(3\ts 2)}(a)\ts
c_{(2)\tss(1)}^{\tss\rho}(a)
=c_{(3)\tss(2)}^{\tss(3\ts 2)}(a)\ts c_{(2)\tss(1)}^{\tss(2)}(a)
+c_{(3)\tss(2\ts 1)}^{\tss(3\ts 2)}(a)
\ts c_{(2)\tss(1)}^{\tss(2\ts 1)}(a).
\een
Now, $c_{(3)\tss(2)}^{\tss(3\ts 2)}(a)=c_{(3)\tss(2)}^{\tss(3\ts 2)}=1$
and $c_{(2)\tss(1)}^{\tss(2\ts 1)}(a)=c_{(2)\tss(1)}^{\tss(2\ts 1)}=1$,
while applying \cite[Theorem~2.1]{m:lr} we get
$c_{(2)\tss(1)}^{\tss(2)}(a)=a_{-1}-a_1$ and
$c_{(3)\tss(2\ts 1)}^{\tss(3\ts 2)}(a)=a_{-2}-a_2$. Hence,
\ben
\wh K_{(3\ts 2)\tss(3\ts 2\ts 1)}(a)=
a_{-2}+a_{-1}-a_1-a_2.
\een
\vskip-1.2\baselineskip
\qed
\eex

The polynomials $K_{\la\mu}(a)$ can be calculated
by the following procedure. Given a partition
$\mu=(\mu_1,\dots,\mu_l)$,
write each dual complete
symmetric function $\wh h_{\mu_i}(y\vt a)$ as a series of the hook
Schur functions with coefficients in $\QQ[a]$ using
Corollary~\ref{cor:sone}. Then multiply the Schur
functions using the classical Littlewood--Richardson rule.
Finally, use Theorem~\ref{thm:invtabre} to
represent each Schur function as a series of the
dual Schur functions.

\bex\label{ex:konadu}
By Example~\ref{ex:sone}, $\wh h_1(y\vt a)^2$ equals
\ben
\big(s_{(1)}(y)+a_0\tss s_{(2)}(y)-
a_1\tss s_{(1^2)}(y)+a^2_0\tss s_{(3)}(y)
-a_0\tss a_1\tss s_{(2\ts 1)}(y)+a^2_1\tss s_{(1^3)}(y)+\dots\big)^2.
\een
Hence, multiplying the Schur functions, we find that
\ben
\bal
\wh h_1(y\vt a)^2&=s_{(2)}(y)+s_{(1^2)}(y)+2\tss a_0\tss s_{(3)}(y)
+2\tss(a_0-a_1)\tss s_{(2\ts 1)}(y)-2\tss a_1\tss s_{(1^3)}(y)\\
{}&+3\tss a^2_0\tss s_{(4)}(y)
+(3\tss a^2_0-4\tss a_0\tss a_1)\tss s_{(3\ts 1)}(y)
+(a^2_0+a^2_1-2\tss a_0\tss a_1)\tss s_{(2^2)}(y)\\
{}&+(3\tss a^2_1-4\tss a_0\tss a_1)\tss s_{(2\ts 1^2)}(y)
+3\tss a^2_1\tss s_{(1^4)}(y)+\cdots.
\eal
\een
Expanding now each Schur function with the use of
Theorem~\ref{thm:invtabre} or Corollary~\ref{cor:inext},
we come to
\ben
\bal
\wh h_1(y\vt a)^2&=\wh s_{(2)}(y\vt a)+\wh s_{(1^2)}(y\vt a)
+(a_0-a_{-1})\tss \wh s_{(3)}(y\vt a)
+(a_0-a_1)\tss \wh s_{(2\ts 1)}(y\vt a)\\
{}&+(a_2-a_1)\tss \wh s_{(1^3)}(y\vt a)
+(a_0-a_{-2})\tss(a_0-a_{-1})\tss \wh s_{(4)}(y\vt a)\\
{}&+(a_0-a_1)\tss(a_0-a_{-1})\tss \wh s_{(3\ts 1)}(y\vt a)
+(a_0-a_1)^2\tss  \wh s_{(2^2)}(y\vt a)\\
{}&+(a_1-a_0)\tss(a_1-a_2)\tss \wh s_{(2\ts 1^2)}(y\vt a)
+(a_1-a_2)\tss(a_1-a_3)\tss \wh s_{(1^4)}(y\vt a)+\cdots,
\eal
\een
thus calculating the first few polynomials $K_{\la\tss(1^2)}(a)$.
\qed
\eex

\bex\label{ex:monom}
Using Example~\ref{ex:konadu}, we can calculate the first few
double monomial symmetric functions:
\ben
\bal
m_{(1)}(x\vt a)&=s_{(1)}(x\vt a),\qquad
m_{(1^2)}(x\vt a)=s_{(1^2)}(x\vt a)\\
m_{(2)}(x\vt a)&=s_{(2)}(x\vt a)-s_{(1^2)}(x\vt a)\\
m_{(1^3)}(x\vt a)&=s_{(1^3)}(x\vt a)
+(a_1-a_2)\tss s_{(1^2)}(x\vt a)\\
m_{(2\ts 1)}(x\vt a)&=s_{(2\ts 1)}(x\vt a)-2\tss s_{(1^3)}(x\vt a)
+(2\tss a_2-a_1-a_0)\tss s_{(1^2)}(x\vt a)\\
m_{(3)}(x\vt a)&=s_{(3)}(x\vt a)-
s_{(2\ts 1)}(x\vt a)+\tss s_{(1^3)}(x\vt a)
+(a_{-1}-a_0)\tss s_{(1^2)}(x\vt a).
\eal
\een
\vskip-1.1\baselineskip
\qed
\eex

The following formula for the dual character
polynomials is implied by
Corollary~\ref{cor:expasdouble}.

\bco\label{cor:char}
We have
\ben
\wh\chi^{\ts\la}_{\mu}(a)=\sum_{\rho}
(-1)^{m(\la/\rho)}\ts \chi^{\rho}_{\mu}\ts
\psi_{\la/\rho}(a),
\een
summed over diagrams $\rho$ with $|\rho|=|\mu|$.
\qed
\eco

\section{Interpolation formulas}
\label{sec:idd}
\setcounter{equation}{0}

\subsection{Rational expressions for the transition coefficients}

Applying Proposition~\ref{prop:interp}, we can get
expressions for the polynomials
$\wh K_{\la\mu}(a)$, $\chi_{\mu}^{\la}(a)$,
$\wh c_{\la\mu}^{\ts\tss\nu}(a)$
and $c_{\la\mu}^{\tss\nu}(a)$ as rational
functions in the variables $a_i$.

\bpr\label{prop:ratfe}
We have the expressions
\begin{align}
\label{kostka}
\wh K_{\la\mu}(a)&=\sum_{R}\sum_{k=0}^{l}
\frac{h_{\mu}(a_{\rho^{(k)}}\vt a)}{(|a_{\rho^{(k)}}|-|a_{\rho^{(0)}}|)
\ldots\wedge\ldots(|a_{\rho^{(k)}}|-|a_{\rho^{(l)}}|)},\\
\label{charact}
\chi_{\mu}^{\la}(a)&=\sum_{R}\sum_{k=0}^{l}
\frac{p_{\mu}(a_{\rho^{(k)}}\vt a)}{(|a_{\rho^{(k)}}|-|a_{\rho^{(0)}}|)
\ldots\wedge\ldots(|a_{\rho^{(k)}}|-|a_{\rho^{(l)}}|)},\\
\label{duallr}
\wh c_{\la\mu}^{\ts\tss\nu}(a)&=\sum_{R}\sum_{k=0}^{l}
\frac{s_{\nu/\mu}(a_{\rho^{(k)}}\vt a)}{(|a_{\rho^{(k)}}|-|a_{\rho^{(0)}}|)
\ldots\wedge\ldots(|a_{\rho^{(k)}}|-|a_{\rho^{(l)}}|)},
\end{align}
summed over all sequences of partitions $R$
of the form
\ben
\varnothing=\rho^{(0)}\to\rho^{(1)}\to
\dots\to\rho^{(l-1)}\to\rho^{(l)}=\la.
\een
Moreover,
\beql{lrpo}
c_{\la\mu}^{\tss\nu}(a)=\sum_{R}\sum_{k=0}^{l}
\frac{s_{\la}(a_{\rho^{(k)}}\vt a)}{(|a_{\rho^{(k)}}|-|a_{\rho^{(0)}}|)
\ldots\wedge\ldots(|a_{\rho^{(k)}}|-|a_{\rho^{(l)}}|)},
\eeq
summed over all sequences of partitions $R$
of the form
\ben
\mu=\rho^{(0)}\to\rho^{(1)}\to
\dots\to\rho^{(l-1)}\to\rho^{(l)}=\nu.
\een
\vskip-1.2\baselineskip
\qed
\epr

The last formula was given in \cite{ms:lr}
for polynomials closely related to $c_{\la\mu}^{\tss\nu}(a)$.
Due to \eqref{cla}, the Kostka numbers $K_{\la\mu}$
and the values of the irreducible
characters $\chi^{\la}_{\mu}$ of the symmetric group
can be found from \eqref{kostka} and \eqref{charact}.

\bex\label{ex:valide}
If $|\la|=n$, then
\ben
\chi^{\la}_{(1^n)}=\sum_{R}\sum_{k=1}^{n}
\frac{(|a_{\rho^{(k)}}|-|a_{\rho^{(0)}}|)^{n-1}}
{(|a_{\rho^{(k)}}|-|a_{\rho^{(1)}}|)
\ldots\wedge\ldots(|a_{\rho^{(k)}}|-|a_{\rho^{(n)}}|)}
=\sum_{R}\ts 1,
\een
which coincides with the
number of standard
$\la$-tableaux.
\qed
\eex

\subsection{Identities with dimensions of skew diagrams}

Specializing the variables by setting $a_i=-i+1$ for all $i\in\ZZ$
in the expressions of Proposition~\ref{prop:ratfe},
we obtain some identities for
the Kostka numbers,
the values of the irreducible
characters and the Littlewood--Richardson coefficients
involving dimensions of skew diagrams.
Under this specialization, the double symmetric functions
become the {\it shifted symmetric functions\/} of \cite{oo:ss},
so that some of the combinatorial
results concerning the ring $\La(x\vt a)$
discussed above in the paper reduce to
the respective results of
\cite{oo:ss} for the ring $\La^*$
of shifted symmetric functions; see also \cite{io:kc}
for an alternative description of the ring $\La^*$.

For any
skew diagram $\theta$ denote by
$\dim\theta$
the number of standard $\theta$-tableaux
(i.e., row and column strict) with
entries in $\{1,2,\dots,|\theta|\}$ and set
\ben
H_{\theta}=\frac{|\theta|!}{\dim\theta}.
\een
If $\theta$ is normal (nonskew), then
$H_{\theta}$ coincides with the product of
the hooks of $\theta$ due to the hook formula.
Under the specialization $a_i=-i+1$, for any partition
$\mu$ we have
\ben
a_{\mu}=(\mu_1,\mu_2-1,\dots).
\een
The following formula for the values
of the double Schur functions was proved in \cite{oo:ss}:
if $\mu\subseteq \nu$, then
\ben
s_{\mu}(a_{\nu}\vt a)=\frac{H_{\nu}}{H_{\nu/\mu}}.
\een
This formula is deduced from Proposition~\ref{prop:interp}
with the use of \eqref{hoo} which takes the form
$s_{\la}(a_{\la}\vt a)=H_{\la}$. Then \eqref{lrpo}
implies the identity for the
Littlewood--Richardson coefficients $c_{\la\mu}^{\tss\nu}$
which was proved in \cite{ms:lr}:
\ben
c_{\la\mu}^{\nu}=\sum_{\rho}(-1)^{|\nu/\rho|}
\frac{H_{\rho}}{H_{\nu/\rho}\ts H_{\rho/\la}\ts H_{\rho/\mu}},
\een
summed over diagrams $\rho$ which contain both $\la$ and $\mu$,
and are contained in $\nu$.

We also have the respective consequences of
\eqref{kostka} and \eqref{charact}.
For partitions
$\mu=(1^{m_1}2^{m_2}\dots r^{m_r})$ and
$\rho=(\rho_1,\dots,\rho_l)$ set
\ben
\pi_{\mu}(\rho)=\prod_{k=1}^r
\Big((1-\rho_1)^{k}+\dots+
(l-\rho_l)^{k}-1^k-\dots-l^{\tss k}\Big)^{m_k}
\een
and
\ben
\varkappa_{\mu}(\rho)=\prod_{k=1}^r
\Big(\sum_{i_1\geqslant\dots\geqslant i_k\geqslant 1}
\rho_{i_1}(\rho_{i_2}-1)\dots (\rho_{i_k}-k+1)\Big)^{m_k}.
\een

The following formulas are obtained by specializing
$a_i=i$ and $a_i=-i+1$, respectively, in \eqref{charact}
and \eqref{kostka}.

\bco\label{cor:kochar}
Let $\la$ and $\mu$ be partitions of $n$.
Then
\ben
\chi^{\la}_{\mu}=\sum_{\rho\subseteq\la}\ts
\frac{(-1)^{|\rho|}\ts \pi_{\mu}(\rho)}{H_{\rho}\ts H_{\la/\rho}}
\een
and
\ben
K_{\la\mu}=\sum_{\rho\subseteq\la}\ts
\frac{(-1)^{|\la/\rho|}\ts
\varkappa_{\mu}(\rho)}{H_{\rho}\ts H_{\la/\rho}}.
\een
\vskip-1.2\baselineskip
\qed
\eco

\bex\label{ex:char}
Let
$\la=(3\ts 2)$ and $\mu=(2\ts 1^3)$. Then
\ben
\pi_{\mu}(\rho)=-(\rho_1+\rho_2)^3\ts
(\rho_1^{\tss 2}+\rho_2^{\tss 2}-2\ts\rho^{}_1-4\ts\rho^{}_2),
\een
and
\ben
\bal
H_{(3\ts 2)/(1)}&=24/5,\qquad
H_{(3\ts 2)/(2)}=2,\qquad
H_{(3\ts 2)/(1^2)}=3,\qquad
H_{(3\ts 2)/(3)}=2,\\
H_{(3\ts 2)/(2\ts 1)}&=1,\qquad
H_{(3\ts 2)/(2\ts 2)}=1,\qquad H_{(3\ts 2)/(3\ts 2)}=1.
\eal
\een
Hence,
\ben
\chi^{(3\ts 2)}_{(2\ts 1^3)}=
-\frac{5}{24}+\frac{32}{6}+\frac{81}{12}
-\frac{81}{3}+\frac{256}{12}-\frac{125}{24}=1.
\een
\eex

\end{document}